\pdfoutput=1

\newif\ifreview
\reviewtrue
\ifreview
\documentclass[10pt]{article}
\else
\documentclass[10pt,twocolumn]{article}
\fi
\usepackage{lmodern}
\usepackage{microtype}

\usepackage{cite}
  \usepackage[pdftex]{graphicx}

\usepackage{savesym}
\savesymbol{iint}
\savesymbol{iiint}

\usepackage[cmex10]{amsmath}
\usepackage{amsfonts}

\restoresymbol{AMS}{iint}
\restoresymbol{AMS}{iiint}
\usepackage{algorithmic}
\usepackage{array}
  \usepackage[caption=false,font=footnotesize]{subfig}
\usepackage{fixltx2e}

\usepackage{bm}

\usepackage{hyperref}
\usepackage[capitalise]{cleveref}
\usepackage{booktabs}
\usepackage{algorithm}
\usepackage{breqn}

\begin{document}
\ifreview
\newcommand{\figwidth}{0.32\linewidth}
\newcommand{\figwidthfour}{0.24\linewidth}
\else
\newcommand{\figwidth}{0.32\linewidth}
\newcommand{\figwidthfour}{0.24\linewidth}
\fi
\title{A method for locally approximating regularized iterative tomographic reconstruction methods}

\author{Dani\"{e}l~M.~Pelt*
        and~Kees~Joost~Batenburg%
\thanks{D.M. Pelt is currently a PhD candidate at CWI, P.O. Box 94079, 1090 GB, Amsterdam. (email: D.M.Pelt@cwi.nl tel:+31205924162 fax:+31205924199.)}%
\thanks{K.J. Batenburg is with CWI, Amsterdam and with the Mathematical Institute, Leiden University, The Netherlands. (email: joost.batenburg@cwi.nl.)}%
}

\maketitle

\begin{abstract}
    In many applications of tomography, the acquired projections are either limited in number or contain a significant amount of noise.
    In these cases, standard reconstruction methods tend to produce artifacts that can make further analysis difficult.
    Advanced regularized iterative methods, such as total variation minimization, are often able to achieve a higher reconstruction quality by exploiting prior knowledge about the
    scanned object.
    In practice, however, these methods often have prohibitively long computation times or
    large memory requirements.
    Furthermore, since they are based on minimizing a global objective function, regularized iterative methods need to reconstruct the entire scanned object, even when one
    is only interested in a (small) region of the reconstructed image.

    In this paper, we present a method to approximate regularized iterative reconstruction methods inside
    a (small) region of the scanned object. The method only performs computations inside the region of interest, ensuring
    low computational requirements.
    Reconstruction results for different phantom images and types of regularization are given,
    showing that reconstructions of the proposed local method are almost identical to those of the global regularized
    iterative methods that are approximated, even for relatively small regions of interest. Furthermore, we show that larger regions
    can be reconstructed efficiently by reconstructing several small regions in parallel and combining them into a single reconstruction afterwards.
\end{abstract}

\section{Introduction}
The goal of tomography is to
reconstruct an object given its projections for different angles.
Using tomography, it is possible to nondestructively examine
the interior of objects, which makes it useful for many applications.
Examples of tomography in practice include computed tomography in medicine and
electron tomography in materials science.
Because of its
practical usefulness, many algorithms have been developed to perform tomographic
reconstruction. An overview of past research on tomography can be found in
\cite{Kak2001,Natterer2001,Buzug2008}.
Two types of reconstruction methods are commonly used: \emph{analytical} methods,
which discretize a continuous inversion formula of the problem, and \emph{algebraic}
methods, in which a linear system that represents the problem is solved.

In many applications of tomography, it is impossible to acquire a large number
of low-noise projections. For example, when scanning live animals, there is a
limit on the total dose deposited on the animal during the experiment
\cite{Lovric2013}. In electron tomography, the scanned sample is damaged by the
electron beam, which leads to a limit on the number of projections that can be
acquired \cite{McEwen1995}.
In these cases, standard reconstruction methods often fail to
produce reconstructions with adequate quality for further analysis
\cite{Lovric2013}.
Analytical methods are based on the assumption that noise-free projections are available
for \emph{all} angles, which is infeasible in practice.
In algebraic methods,
the linear system that is solved is typically both underdetermined and ill-conditioned, which can make it
difficult to find accurate reconstructions when the available projection data is limited
and/or noisy.

Recently developed advanced reconstruction methods aim to improve reconstruction
quality by exploiting prior knowledge about the scanned object or scanning
system. Often, these methods add additional terms to the objective function that
is minimized in standard algebraic methods.
Methods of this type will be called \emph{regularized iterative methods} in this paper.
For example, if it is known beforehand that the physical quantity that is
reconstructed cannot be negative, a nonnegativity constraint can be added
to the objective function to improve the reconstruction quality.
If it is known that the scanned object has a sparse boundary, total variation
minimization can be applied by adding a term that minimizes the gradient
of the reconstructed image~\cite{Sidky2008}.
If the added prior knowledge is appropriate for the acquired
data, regularized iterative methods can be extremely successful in reconstructing objects
from (highly) limited data \cite{Batenburg2011,Kostenko2013}.

One of the main disadvantages of regularized iterative methods is their computational
cost, which is typically very high. A high computational cost of a
reconstruction method can be prohibitive for its application in practice. For
example, in ultrafast tomographic experiments at synchrotrons, the computation time
of the reconstruction method has to match the high speed of the acquisition
of projection data \cite{Mokso2013}.
An additional problem is that regularized iterative methods often have a number of tunable
parameters that influence the reconstruction quality greatly.
In many cases, values for these parameters are
chosen by trial-and-error, which can be very time-consuming for methods with a
high computational cost.
These problems are especially important in cases where a large object is scanned,
but the features of interest are only located in a small region of the object.
Since
regularized iterative methods, and the algebraic methods
they are based on, minimize a \emph{global} objective function, they
typically need to compute the entire volume
during reconstruction, which may not fit into the available memory of
the graphic processing units used
to perform the reconstruction \cite{Xu2005}.

Analytical methods, on the other
hand, can be evaluated \emph{locally}: if one is only interested in a small
subvolume of the reconstruction, only that subvolume has to be reconstructed.
When reconstructing large volumes, analytical methods can divide the
reconstruction volume into subvolumes that do fit into the available memory, and
reconstruct each subvolume separately, resulting in an efficient method to
compute the full reconstruction volume.
This property is one of the reasons that in many applications of tomography,
standard analytical methods are still the most popular reconstruction methods
instead of regularized iterative methods \cite{Pan2009}.

In previous research, attempts have been made to improve the computation time of
algebraic methods when one is only interested in a small region. One approach is based
on subtracting simulated projections of a global reconstruction outside the
region of interest from the acquired projections~\cite{Ziegler2008,DeWitte2010,Kopp2015}. The resulting
altered projections are used in an algebraic reconstruction of the region of interest.
If the quality of the initial global reconstruction is not sufficient, however, this approach
can result in artifacts inside the reconstructed region. Furthermore, the global analytical reconstruction
has to fit into the available memory of the computational system, which may pose problems for large
objects. Another approach is to use a multiresolution algebraic method~\cite{Niinimaeki2007},
computing a global reconstruction with a higher resolution inside the region of interest and a lower resolution outside the region.
One problem with this approach is that the resolution outside the region of interest has to
be sufficiently high to prevent artifacts inside the region of interest, which may lead to
large computational costs. Choosing the correct resolution parameters to use can be difficult,
since it depends on the scanned object and the acquisition parameters.
Note that both approaches are typically used for approximating algebraic methods without any additional regularization terms.

In this paper, we present a novel approach, resulting in a method for approximating a
computationally expensive regularized iterative method in a (small) subvolume of the full
reconstruction volume.
The proposed method only performs computations in the chosen subvolume, ensuring
low computational and memory requirements.
If one is only interested in part of the scanned object,
the new method can significantly reduce the time needed to reconstruct that part
compared to existing regularized iterative methods.
If one wants to reconstruct the entire object, the proposed
method also allows for significant reduction of computation time by enabling
parallel computation of different subvolumes, and it enables regularized iterative
reconstruction of large datasets that
do not fit completely into the available memory.
In addition, the method can be
used to quickly estimate parameters of a slow regularized iterative method by estimating
them in a small subvolume.

The proposed method is based on approximating standard algebraic methods by
a modified analytical method. In recent years, several methods have been
proposed that achieve this by modifying the filter that is typically used
in analytical methods.
In one study, an angle-independent filter is calculated based on analytic
analysis of the algebraic SIRT method \cite{Zeng2012}. An extension of the
method for noisy projection data is given in \cite{Zeng2013}. In another study,
a method of calculating a data-dependent filter is given in \cite{Pelt2014}.
Finally, an angle-dependent and geometry-dependent filter is calculated by
repeated application of the SIRT method in \cite{Batenburg2012}. A faster method of
calculating similar filters for the algebraic SIRT method is proposed in
\cite{Pelt2015}.
None of these methods, however, allow for
inclusion of popular prior knowledge terms, such as total variation
minimization, which can limit their usefulness in practice.

We first show
the application of the filter of~\cite{Pelt2015} to locally approximate the algebraic SIRT
method.
Then, we extend the method to allow for
local approximation of a regularized iterative method as well.
Finally, we demonstrate that the proposed method is able to produce local reconstructions that
are very similar to reconstructions of global regularized iterative methods for various
types of exploited prior knowledge.

This paper is structured as follows. In \cref{sec:notconc}, we introduce the
notations we use throughout the paper, and formally define the tomographic
reconstruction problem and the standard analytical and algebraic approaches. The
main contribution of this paper is given in \cref{sec:meth}, where we first
apply the method proposed in \cite{Pelt2015} to approximate SIRT locally. We
then extend this approximation by including prior knowledge in the
reconstruction of a subvolume, and give some details on how to implement
the resulting method in practice.
The experiments we performed to study the new
method are explained in \cref{sec:exp}, and the results of those experiments are
shown in \cref{sec:res}. We conclude in \cref{sec:conc} with a brief summary of
the paper and some final remarks.

\section{Notation and concepts}
\label{sec:notconc}
In this section, the mathematical notation that we use throughout the paper is
introduced, and a formal definition of the tomographic reconstruction problem is
given. The standard analytical and algebraic approaches to the problem are
explained, and their mathematical definitions are given. Finally, we explain how
prior knowledge can be exploited in algebraic methods by extending their
objective functions, resulting in regularized iterative methods.

\subsection{Notation and problem definition}
\begin{figure}
\begin{center}
\includegraphics{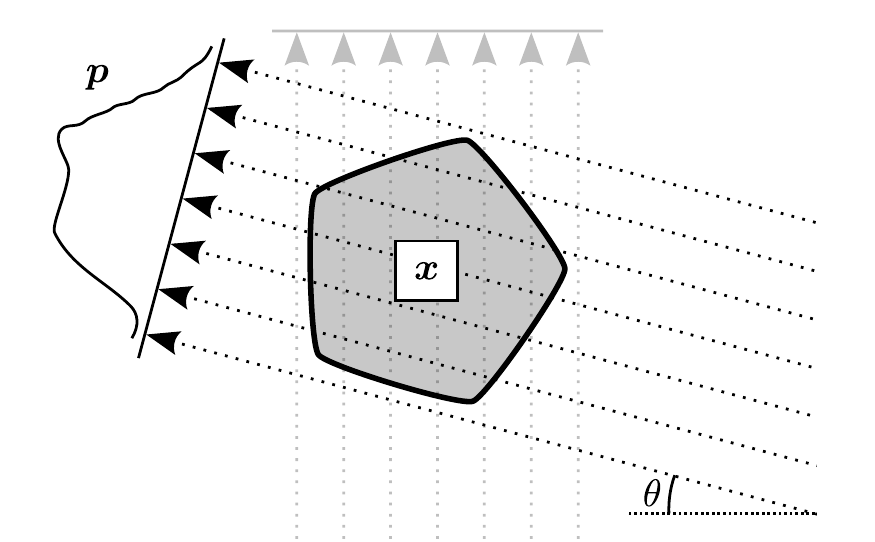}
\end{center}
\caption{The two-dimensional parallel-beam geometry used in this paper. Parallel lines, rotated by angle $\theta$, pass through the object $\boldsymbol{x}$. The projections $\boldsymbol{p}$ of $\boldsymbol{x}$ are given by the line integrals of $\boldsymbol{x}$ over the lines.}
\label{fig:parbeam}
\end{figure}

We focus on two-dimensional parallel-beam tomographic reconstruction problems with a single rotation axis.
Note that in many cases it is possible to convert other tomographic geometries, such as cone-beam or spiral tomography,
to a parallel-beam geometry by rebinning~\cite{Grass2000,Kachelries2000}.
Parallel-beam projection data are acquired by rotating an array of detectors around the object (or, equivalently, rotating the object), with the detectors of the array located on a straight line. This acquisition scheme is shown graphically in \cref{fig:parbeam}. If the number of detectors in the array is denoted by $N_d$, and the number of rotation angles for which data are acquired is denoted by $N_\theta$, we can write the measured line integrals as a vector $\boldsymbol{p}$ with $N_d N_\theta$ elements, one for each combination of detector element and rotation angle. The reconstructed image is represented as a vector $\boldsymbol{x}$ with $N^2$ elements, one for each pixel of the $N \times N$ pixel grid on which the reconstruction is calculated. The main problem in tomographic reconstruction is to find the unknown image $\boldsymbol{x}$, given the acquired projection data $\boldsymbol{p}$.

The \emph{forward projection} operator $\mathcal{W}: \mathbb{R}^{N^2} \rightarrow \mathbb{R}^{N_d N_{\theta}}$ is the operator that, for a given projection geometry, corresponds to the discretized line integrals of an object represented on a $N \times N$ pixel grid. Using the above notation, we can write this operator as a $N_d N_{\theta} \times N^2$ matrix $\boldsymbol{W}$, with element $w_{ij}$ giving the contribution of pixel $j$ to detector $i$. The transpose of this operator, $\boldsymbol{W}^T$, is called the \emph{backprojection} operator. Typically, a forward projection of an image $\boldsymbol{x}$ is calculated on-the-fly by calculating its line integrals directly \cite{Palenstijn2011}. Similarly, multiplying $\boldsymbol{p}$ by $\boldsymbol{W}^T$ is done implicitly by backprojecting $\boldsymbol{p}$ on-the-fly. The advantage of this approach is that the matrix $\boldsymbol{W}$, which can be very large, never has to be stored in memory. Furthermore, forward projections and backprojections can be computed very efficiently on graphic processor units (GPUs) \cite{Xu2005,Mueller2007}.

Our novel approach aims to reconstruct only a \emph{local} part $\mathcal{L}$ of the entire reconstruction grid. Here, $\mathcal{L}$ is a subset of all $N^2$ pixels of the entire reconstruction grid, usually ordered in a $N_{\mathcal{L}}\times N_{\mathcal{L}}$ grid as well. Let $\boldsymbol{M}_{\mathcal{L}}$ be a diagonal matrix with a value $1$ on the diagonal of row $i$ if pixel $i$ is inside ${\mathcal{L}}$, and $0$ everywhere else. In other words, $\boldsymbol{M}_{\mathcal{L}}$ keeps all pixels of an image that are inside ${\mathcal{L}}$, and zeros all other pixels. Similarly, we define a matrix
$\boldsymbol{M}_{\mathcal{F}}$ that zeros all pixels inside ${\mathcal{L}}$, and keeps all other pixels. Using these, we can define local operators $\boldsymbol{W}_{\mathcal{L}}$ and $\boldsymbol{W}_{\mathcal{L}}^T$, and outer operators $\boldsymbol{W}_{\mathcal{F}}$ and $\boldsymbol{W}_{\mathcal{F}}^T$:
\begin{equation}
  \begin{aligned}
\boldsymbol{W}_{\mathcal{L}} &= \boldsymbol{W} \boldsymbol{M}_{\mathcal{L}} \\
\boldsymbol{W}_{\mathcal{L}}^T &= \boldsymbol{M}_{\mathcal{L}} \boldsymbol{W}^T \\
\boldsymbol{W}_{\mathcal{F}} &= \boldsymbol{W} \boldsymbol{M}_{\mathcal{F}} \\
\boldsymbol{W}_{\mathcal{F}}^T &= \boldsymbol{M}_{\mathcal{F}} \boldsymbol{W}^T
  \end{aligned}
\end{equation}
Since $\boldsymbol{M}_{\mathcal{L}} + \boldsymbol{M}_{\mathcal{F}} = \boldsymbol{I}$ by construction, we have that the sum of $\boldsymbol{W}_{\mathcal{L}}$ and $\boldsymbol{W}_{\mathcal{F}}$ is equal to $\boldsymbol{W}$:
\begin{equation}
  \boldsymbol{W} = \boldsymbol{W}_{\mathcal{L}} + \boldsymbol{W}_{\mathcal{F}}
  \label{eq:locouter}
\end{equation}
Note that local forward projections and backprojections can be computed significantly faster than full forward projections and backprojections, since many rows and columns of $\boldsymbol{W}_{\mathcal{L}}$ and $\boldsymbol{W}_{\mathcal{L}}^T$ are zero.

\subsection{Common reconstruction methods}
\label{sec:introrecs}

Using the above definitions, we can write one of the most popular reconstruction
methods, the analytical filtered backprojection (FBP) method, as:
\begin{equation}
  \mathit{FBP}(\boldsymbol{p},\boldsymbol{h}) = \boldsymbol{W}^T \boldsymbol{C}_{\boldsymbol{h}} \boldsymbol{p} \label{eq:fbp}
\end{equation}
Here, $\boldsymbol{C}_{\boldsymbol{h}}$ is a convolution operator that
convolves each 1D array of detector values, taken at a single rotation angle,
with the 1D filter $\boldsymbol{h}$~\cite{Kak2001}. Note that this 1D filter can be
different for each rotation angle. Several fixed angle-independent filters are commonly used in
practice, such as the Ram-Lak (ramp), Shepp-Logan, and Hann filters
\cite{Farquhar1997}. One reason for the popularity of FBP is its computational
efficiency: the filtering step can be performed very efficiently in Fourier
space, and only one backprojection has to be computed during reconstruction.
Another advantage of the filtered backprojection method compared to other
methods is that we can calculate its values inside the local part
${\mathcal{L}}$ by simply exchanging $\boldsymbol{W}^T$ by
$\boldsymbol{W}_{\mathcal{L}}^T$ in \cref{eq:fbp}:
\begin{equation}
  \mathit{FBP}_{\mathcal{L}}(\boldsymbol{p},\boldsymbol{h}) = \boldsymbol{W}_{\mathcal{L}}^T \boldsymbol{C}_{\boldsymbol{h}} \boldsymbol{p} \label{eq:fbplocal}
\end{equation}

A different approach to solving the reconstruction problem is the algebraic
approach. Here, we form a linear system $\boldsymbol{W} \boldsymbol{x} =
\boldsymbol{p}$, and solve for $\boldsymbol{x}$. Most algebraic methods find a
solution $\boldsymbol{x}_{\text{alg}}$ by minimizing the difference, in some
vector norm, between the forward projection of the solution and the measured
projection data. This difference is called the \emph{projection error}. In the
case of the $\ell_2$-norm, we can write this as:
\begin{equation}
\boldsymbol{x}_{\mathit{alg}} = \underset{\boldsymbol{x}}{\operatorname{argmin}} \left\lVert \boldsymbol{p} - \bm{W} \boldsymbol{x}\right\rVert_2^2
\label{eq:algsys}
\end{equation}
Since the matrix $\boldsymbol{W}$ is often very large, \cref{eq:algsys} is
usually not solved directly. Instead, an iterative optimization method is
typically used to iteratively decrease the projection error. Implicit
regularization of the solution can be included by stopping the iteration process
early, which is needed because $\boldsymbol{W}$ is usually ill-conditioned and
noise is often present in $\boldsymbol{p}$.

Different iterative optimization methods can be used to minimize the projection
error, leading to different algebraic methods. The CGLS method, for example, is
based on a conjugate gradient method \cite{Bjoerck1996}. Another popular
algebraic method is the simultaneous iterative reconstruction technique (SIRT)
\cite{Kak2001}. The SIRT method belongs to the class of Landweber iteration
methods \cite{Landweber1951}, and uses a specific Krylov subspace method to
minimize the projection error iteratively. A single iteration of the SIRT method
can be viewed as a gradient-descent step on the projection error, and can be written as:
\begin{equation}
  \boldsymbol{x}_s^{k+1} = S(\boldsymbol{x}_s^k) = \boldsymbol{x}_s^k + \alpha \boldsymbol{W}^T \left( \boldsymbol{p} - \boldsymbol{W} \boldsymbol{x}_s^k \right) \label{eq:sirtit}
\end{equation}
Note that in algebraic methods, we are not able to simply exchange
$\boldsymbol{W}$ by $\boldsymbol{W}_{\mathcal{L}}$ to find the reconstruction
inside ${\mathcal{L}}$, since then we would be solving the linear system
$\boldsymbol{W}_{\mathcal{L}} \boldsymbol{x} = \boldsymbol{p}$, which will have
a completely different solution than $\boldsymbol{W} \boldsymbol{x} =
\boldsymbol{p}$ if the scanned object is nonzero outside $\mathcal{L}$.

\subsection{Regularized iterative methods}
\begin{figure}
\begin{center}
\subfloat[][]{\includegraphics[width=\figwidthfour]{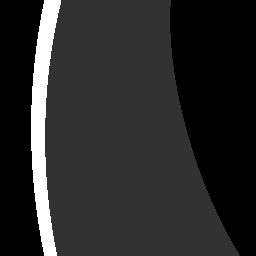}}
\hfill
\subfloat[][]{\includegraphics[width=\figwidthfour]{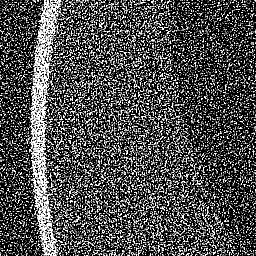}}
\hfill
\subfloat[][]{\includegraphics[width=\figwidthfour]{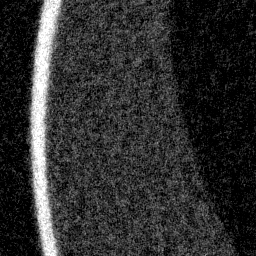}}
\hfill
\subfloat[][]{\includegraphics[width=\figwidthfour]{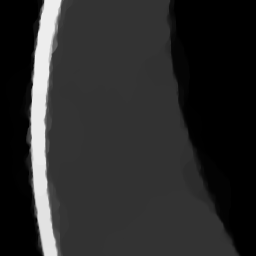}}
\end{center}
\caption{Zoomed-in reconstructions of the Shepp-Logan head phantom (a), showing the resulting images of three different reconstruction methods: (b) FBP, (c) SIRT, and (d) total variation minimization. The images were reconstructed on a $1024 \times 1024$ pixel grid, using projection data acquired with $N_d=1024$ detectors and $N_\theta=256$ projection angles, equally distributed in the interval $[0,\pi]$, and additional Poisson noise applied.}
\label{fig:comp}
\end{figure}

A common way of including prior knowledge in algebraic methods is to add
additional constraints to the objective function of \cref{eq:algsys}. In this
paper, we distinguish two types of constraints that are commonly used:
\emph{domain constraints}, which restrict the domain of possible solutions,
and \emph{penalty constraints}, which penalize undesired solutions in the
objective function. The resulting regularized iterative reconstructions can be written as:
\begin{equation}
\boldsymbol{x}_{\mathit{reg}} = \underset{\boldsymbol{x} \in D}{\operatorname{argmin}} \left[ \left\lVert \boldsymbol{p} - \bm{W} \boldsymbol{x}\right\rVert_2^2 + \lambda g(\boldsymbol{x}) \right]
\label{eq:priorsys}
\end{equation}
Here,
$D$ is a restricted domain for the possible solutions $\boldsymbol{x}$, and
$g: \mathbb{R}^{N^2} \rightarrow \mathbb{R}$ is a penalty function that
penalizes solutions that do not fit with the assumed prior knowledge.
The
$\lambda$ term controls how strongly the penalty function is weighted compared to the
projection error term.
The domain $D$ is used to specify domain constraints, for example when adding
a nonnegativity constraint on the values of $\boldsymbol{x}$ by using $D = \{\boldsymbol{x} \in \mathbb{R}^{N^2}; x_i\geq 0, \, i=1,\dots,N^2\}$.
The cost function $g(\boldsymbol{x})$ is used to specify penalty constraints.
For example, if we assume that the scanned object is sparse in
some wavelet basis, we can set $g(\boldsymbol{x}) = \left\lVert \boldsymbol{B}
\boldsymbol{x} \right\rVert_1$, where $\boldsymbol{B}$ is the wavelet
decomposition operator. Similarly, if we assume that the gradient of the scanned
object is sparse, we set $g(\boldsymbol{x}) = \left\lVert \nabla \boldsymbol{x}
\right\rVert_1$ to perform total variation minimization, where
$\nabla$ is a discrete gradient operator. Several algorithms exist that are able
to find solutions to \cref{eq:priorsys}, such as the popular fast iterative
shrinkage-thresholding algorithm (FISTA) \cite{Beck2009}, Chambolle-Pock algorithms~\cite{Chambolle2011}, and adaptive steepest descent projection onto convex sets algorithm (ASD-POCS)~\cite{Sidky2008}. A comparison of
reconstructions obtained using FBP, SIRT, and total variation minimization from
noisy projection data is shown in \cref{fig:comp}.

Many regularized iterative methods
use a scheme that alternates between
gradient-descent steps
on the projection error $\lVert \boldsymbol{p} - \boldsymbol{W} \boldsymbol{x}
\rVert_2^2$, steps that minimize the penalty function $g(\boldsymbol{x})$,
and steps that enforce the domain constraints $D$.
Since a single iteration of the SIRT method is identical to
a single gradient-descent step on the projection error,
these regularized iterative methods can be viewed as a combination
of SIRT iterations and some additional
steps incorporating the prior knowledge.
As an example, one can include box constraints on the
values of the reconstruction pixels of the form $l \leq x_i \leq r$,
which is a domain constraint with $D = \{ \boldsymbol{x} \in \mathbb{R}^{N^2}; l \leq x_i \leq r, \, i=1,\dots,N^2\}$
by using
the following iterations for pixel $i$ of the reconstruction:
\begin{equation}
  \boldsymbol{x}^{k+1}_i = \left\{
  \begin{array}{lrl}
    l & : & \text{if} \enskip S(\boldsymbol{x}^k)_i < l\\
    r & : & \text{if} \enskip S(\boldsymbol{x}^k)_i > r\\
    S(\boldsymbol{x}^k)_i & : & \text{otherwise}
  \end{array}
  \right.
  \label{eq:sirtbox}
\end{equation}
An example of using a penalty constraint is the ISTA method \cite{Daubechies2004} for $\ell_1$-norm
minimization of a representation of the reconstructed image in a wavelet basis.
In this case, a single iteration of the method can be written as:
\begin{equation}
  \boldsymbol{x}^{k+1} = \boldsymbol{B}^{-1}\mathcal{P}_{\lambda}(\boldsymbol{B}S(\boldsymbol{x}^k))
  \label{eq:sirtista}
\end{equation}
where $\boldsymbol{B}$ is the wavelet decomposition operator, and
$\mathcal{P}_{\lambda}$ the soft thresholding operator with threshold $\lambda$:
\begin{equation}
  \mathcal{P}_{\lambda}(\boldsymbol{y})_i = \left\{
  \begin{array}{lrl}
    \operatorname{sgn}(\boldsymbol{y}_i) (\lvert \boldsymbol{y}_i \rvert - \lambda) & : & \text{if} \enskip \boldsymbol{y}_i > \lambda \\
    0 & : & \text{otherwise}
  \end{array}
  \right.
\end{equation}

In this paper, we propose a method to locally approximate regularized iterative
reconstruction methods that are a combination of SIRT iterations and additional
steps that incorporate the prior knowledge.

\section{Method}
\label{sec:meth}
In this section, we introduce the major contribution of this paper: a local
approximation method for regularized iterative reconstruction methods. We first explain
the method introduced in \cite{Pelt2015} to approximate the algebraic SIRT
method by FBP with a specific geometry-dependent filter, and show how this
approach can be used to approximate SIRT locally as well.
Afterwards, we extend the approximation to include prior knowledge, improving
the reconstruction quality. Finally, we give details on how we implemented
the resulting method for the experiments of \cref{sec:exp}.
\subsection{Local approximation of SIRT}
\label{sec:appsirt}
Recall that a single iteration of the SIRT method can be written as:
\begin{equation}
	\boldsymbol{x}_s^{k+1} = S(\boldsymbol{x}_s^k) = \boldsymbol{x}_s^k + \alpha \boldsymbol{W}^T \left( \boldsymbol{p} - \boldsymbol{W} \boldsymbol{x}_s^k \right) \tag{\ref{eq:sirtit}}
\end{equation}
Here, $\alpha \in \mathbb{R}$ is a parameter that influences the stability and rate of convergence of the method. In the rest of this paper, we use $\alpha = \left(N_\theta N_d\right)^{-1}$.

To find an approximation method for the SIRT method, we start by rewriting the equation of a single SIRT iteration (\cref{eq:sirtit}) in a matrix format:
\begin{equation}
\boldsymbol{x}_s^{k+1} = (\boldsymbol{I} - \alpha \boldsymbol{W}^T \boldsymbol{W}) \boldsymbol{x}_s^k + \alpha \boldsymbol{W}^T \boldsymbol{p}
\label{eq:sirtmat}
\end{equation}
This is a recursion equation of the form $\boldsymbol{x}^{k+1} = \boldsymbol{A} \boldsymbol{x}^k + \boldsymbol{b}$, which has
the following solution for iteration $n$:
\begin{equation}
\boldsymbol{x}_s^{n} = \boldsymbol{A}^n \boldsymbol{x}_s^0 + \alpha \left[ \sum_{k=0}^{n-1} \boldsymbol{A}^k \right] \boldsymbol{W}^T \boldsymbol{p}
\label{eq:sirtfull}
\end{equation}
where $\boldsymbol{A} = \boldsymbol{I} - \alpha \boldsymbol{W}^T\boldsymbol{W}$. Often, the initial image of the SIRT
method is set to the zero image ($\boldsymbol{x}_s^0=\boldsymbol{0}$), in which case we end up with:
\begin{equation}
\boldsymbol{x}_s^{n} = \alpha \left[ \sum_{k=0}^{n-1} \boldsymbol{A}^k \right] \boldsymbol{W}^T \boldsymbol{p}
\label{eq:sirtfull2}
\end{equation}

Now, we want to find a method that can approximate \cref{eq:sirtfull2}. In order to find such a method, we look at the FBP method, and note that, in parallel-beam tomography, convolving a sinogram with a filter and backprojecting the result is identical to backprojecting the sinogram and convolving the resulting image with the backprojected filter:
\begin{equation}
  \mathit{FBP}(\boldsymbol{p},\boldsymbol{h}) = \boldsymbol{H}_{\boldsymbol{h}'} \boldsymbol{W}^T \boldsymbol{p}
  \label{eq:fbp2}
\end{equation}
Here, $\boldsymbol{H}_{\boldsymbol{q}}$ is a 2D convolution with filter $\boldsymbol{q}$, and $\boldsymbol{h} = \boldsymbol{W} \boldsymbol{h}'$.

Note the similarities between the rewritten SIRT equation (\cref{eq:sirtfull2}) and the rewritten FBP equation (\cref{eq:fbp2}), which suggest that we can approximate the SIRT equation by approximating $\sum_{k=0}^{n-1} \boldsymbol{A}^k$ by a 2D convolution operation with filter $\boldsymbol{q}_n$:
\begin{equation}
\boldsymbol{x}_s^n \approx \alpha \boldsymbol{H}_{\boldsymbol{q}_n} \boldsymbol{W}^T \boldsymbol{p}
\end{equation}
A good approximating filter $\boldsymbol{q}_n$ can be found by taking the impulse response of $\sum_{k=0}^{n-1} \boldsymbol{A}^k$:
\begin{equation}
\boldsymbol{q}_n = \sum_{k=0}^{n-1} \boldsymbol{A}^k [0,\dots,0,1,0,\dots,0]^T
\label{eq:filtcalc}
\end{equation}
In other words, we apply $\boldsymbol{A}$ to an image $n-1$ times, starting with an image with only the central pixel set to 1, and sum the resulting images to obtain the 2D filter $\boldsymbol{q}_n$.

Since backprojecting a sinogram and convolving the resulting image is the same as convolving the sinogram with the forward projected filter and backprojecting the result, we can write this as:
\begin{equation}
\begin{aligned}
\boldsymbol{x}_s^n &\approx \boldsymbol{W}^T \boldsymbol{C}_{\boldsymbol{u}_n} \boldsymbol{p} \\
\boldsymbol{u}_n &= \alpha \boldsymbol{W} \boldsymbol{q}_n
\end{aligned}
\label{eq:sirtapp}
\end{equation}
Here, $\boldsymbol{C}_{\boldsymbol{h}}$ is the same convolution operator as in \cref{eq:fbp}, and $\boldsymbol{u}_n$ is the corresponding angle-dependent filter. Comparing \cref{eq:fbp} and \cref{eq:sirtapp}, we conclude that the SIRT method with $n$ iterations can be approximated by the FBP method with a special filter $\boldsymbol{u}_n$:
\begin{equation}
\boldsymbol{x}_s^n \approx \mathit{FBP}(\boldsymbol{p},\boldsymbol{u}_n)
\label{eq:sfbp}
\end{equation}
To summarize, the algorithm to compute an approximating filter is given in \cref{alg:filtcalc}.
For more information on implementing this method, and results for non-local tomographic reconstruction, we refer to \cite{Pelt2015}.

\begin{algorithm}[t]
\begin{algorithmic}
  \REQUIRE $\boldsymbol{W} \in \mathbb{R}^{N_d N_{\theta} \times N^2}$, $n \in \mathbb{Z}^+$, $\alpha \in \mathbb{R}$
  \STATE $\boldsymbol{q}_0  \leftarrow \boldsymbol{0}$
  \STATE $\boldsymbol{c} \leftarrow [0,\dots,0,1,0,\dots,0]^T$
  \FOR{$k=1$ \TO $n$}
  \STATE $\boldsymbol{q}_k  \leftarrow \boldsymbol{q}_{k-1} + \boldsymbol{c}$
  \STATE $\boldsymbol{c}  \leftarrow \boldsymbol{c} - \alpha \boldsymbol{W}^T \boldsymbol{W} \boldsymbol{c}$
  \ENDFOR
  \STATE $\boldsymbol{u}_n \leftarrow \alpha \boldsymbol{W} \boldsymbol{q}_n$
  \RETURN $\boldsymbol{u}_n$
\end{algorithmic}
\caption{Compute an FBP filter that approximates $n$ iterations of SIRT}
\label{alg:filtcalc}
\end{algorithm}

One advantage of this approximation is that, after calculating the filter, the final reconstruction method is identical to standard FBP\@.
Therefore, we can use the same approach as for FBP to evaluate it locally: simply exchanging $\boldsymbol{W}^T$ with $\boldsymbol{W}^T_{\mathcal{L}}$:
\begin{equation}
  \boldsymbol{x}_s^{n} \approx \mathit{FBP}_{\mathcal{L}}(\boldsymbol{p},\boldsymbol{u}_{n})
  \label{eq:locsirtapp}
\end{equation}
Results for locally approximating SIRT with this approach are given in~\cref{sec:locsirtres}.

\subsection{Including regularization}
\label{sec:prior}
As explained in~\cref{sec:introrecs}, many regularized iterative methods include a SIRT
step in their iterative equations. In \cref{sec:appsirt}, we showed that we can
approximate these SIRT steps locally by using the proposed filter method.
However, to locally approximate the complete regularized iterative methods, we need to
perform some extra steps. We start by explicitly splitting the reconstruction
image at iteration $k$ into two parts: a standard SIRT image
$\boldsymbol{x}_s^k$ and a prior-based correction term $\boldsymbol{y}^k$:
\begin{equation}
\boldsymbol{x}^k = \boldsymbol{x}_s^k + \boldsymbol{y}^k
\label{eq:priordef}
\end{equation}
Furthermore, we rewrite the equation for a single iteration of these methods,
such that it consists of a single SIRT step on the previous iteration, and an
additional correction term $\boldsymbol{\mathit{d}}$ that incorporates the
prior knowledge:
\begin{equation}
\boldsymbol{x}^{k+1} = S(\boldsymbol{x}^k) + \boldsymbol{\mathit{d}}^{k+1}
\label{eq:priorit}
\end{equation}
Note that it is usually straightforward to rewrite a regularized iterative method that
uses SIRT to this form, although one would typically not use such a formulation
in practice. For example, SIRT with box constraints (\cref{eq:sirtbox}) can be
written in this form by taking:
\begin{equation}
 \boldsymbol{\mathit{d}}^{k+1}_i = \left\{
   \begin{array}{lrl}
     l-S(\boldsymbol{x}^k)_i & : & \text{if} \enskip S(\boldsymbol{x}^k)_i < l\\
     r-S(\boldsymbol{x}^k)_i & : & \text{if} \enskip S(\boldsymbol{x}^k)_i > r\\
     0 & : & \text{otherwise}
   \end{array}
 \right.
\end{equation}
As another example, iterations of the ISTA method with a wavelet basis
(\cref{eq:sirtista}) can be written in the form of \cref{eq:priorit} by taking:
\begin{equation}
\boldsymbol{\mathit{d}}^{k+1} = \boldsymbol{B}^{-1}\mathcal{P}_{\lambda}(\boldsymbol{B}S(\boldsymbol{x}^k)) - S(\boldsymbol{x}^k)
\end{equation}

Now, we aim to find a local approximation to \cref{eq:priorit}. If we apply a
single SIRT iteration to $\boldsymbol{x}^k$, we get:
\begin{equation}
  \begin{aligned}
S(\boldsymbol{x}^k) &= \boldsymbol{A} \left(\boldsymbol{x}_s^k + \boldsymbol{y}^k \right ) + \alpha \boldsymbol{W}^T \boldsymbol{p} \\
&= \boldsymbol{A} \boldsymbol{x}_s^k + \alpha \boldsymbol{W}^T \boldsymbol{p} + \boldsymbol{A} \boldsymbol{y}^k \\
&= S(\boldsymbol{x}_s^{k}) + \boldsymbol{A} \boldsymbol{y}^k
\end{aligned}
\label{eq:sirtitprior}
\end{equation}
By combining \cref{eq:priorit} and \cref{eq:sirtitprior}, we see that:
\begin{equation}
\boldsymbol{x}^{k+1} = S(\boldsymbol{x}_s^{k}) + \boldsymbol{A} \boldsymbol{y}^k + \boldsymbol{\mathit{d}}^{k+1}
\end{equation}
Using the definition of \cref{eq:priordef}, we can take:
\begin{equation}
  \begin{aligned}
    \boldsymbol{x}_s^{k+1} &= S(\boldsymbol{x}_s^{k}) \\
    \boldsymbol{y}^{k+1} &= \boldsymbol{A} \boldsymbol{y}^k + \boldsymbol{\mathit{d}}^{k+1}
  \end{aligned}
\end{equation}
In order to locally approximate \cref{eq:priorit}, we need to find local
approximations for $\boldsymbol{x}_s^{k+1}$ and $\boldsymbol{y}^{k+1}$.

The iterations of $\boldsymbol{x}_s^{k+1}$ are identical to SIRT iterations, for
which we already derived a local approximation in \cref{sec:appsirt}:
\begin{equation}
  \boldsymbol{x}_s^{k+1} \approx \mathit{FBP}_{\mathcal{L}}(\boldsymbol{p},\boldsymbol{u}_{k+1})
\end{equation}
Furthermore, we can choose to only apply the prior knowledge inside the local
part $\mathcal{L}$. In this case, the prior-based correction term
$\boldsymbol{\mathit{d}}^{k+1}$ is only nonzero for pixels inside
$\mathcal{L}$. To find a local approximation to $\boldsymbol{A}
\boldsymbol{y}^k$, we expand $\boldsymbol{A}$, and use the definition of the
local and outer projection operations \cref{eq:locouter}:
\begin{equation}
  \begin{aligned}
    \boldsymbol{A} \boldsymbol{y}^k &=  \left(\boldsymbol{I} - \alpha \boldsymbol{W}^T \boldsymbol{W} \right) \boldsymbol{y}^k \\
    &= \boldsymbol{y}^k - \alpha \left( \boldsymbol{W}_{\mathcal{L}}^T + \boldsymbol{W}_{\mathcal{F}}^T \right) \boldsymbol{W} \boldsymbol{y}^k \\
    &= \boldsymbol{y}^k - \alpha \boldsymbol{W}_{\mathcal{L}}^T \boldsymbol{W} \boldsymbol{y}^k - \alpha \boldsymbol{W}_{\mathcal{F}}^T \boldsymbol{W} \boldsymbol{y}^k
  \end{aligned}
  \label{eq:ay}
\end{equation}
We approximate \cref{eq:ay} locally by simply ignoring the term $\alpha
\boldsymbol{W}_{\mathcal{F}}^T \boldsymbol{W} \boldsymbol{y}^k$ which affects
the pixels outside $\mathcal{L}$. By ignoring this term, we ignore the effect
that the local prior has on the pixels outside $\mathcal{L}$, which can affect
the pixels inside $\mathcal{L}$ in later iterations. Since we are, in the end,
only interested in the reconstruction inside $\mathcal{L}$, this approximation
is usually sufficiently accurate in practice. Another result of this
approximation is that $\boldsymbol{y}^k$ will be zero outside $\mathcal{L}$ for
any iteration $k$, and therefore we can substitute
$\boldsymbol{W}_{\mathcal{L}}$ for $\boldsymbol{W}$ in the forward projection as well:
\begin{equation}
  \boldsymbol{A} \boldsymbol{y}^k \approx \boldsymbol{y}^k - \alpha \boldsymbol{W}_{\mathcal{L}}^T \boldsymbol{W}_{\mathcal{L}} \boldsymbol{y}^k
\end{equation}

To summarize, we have derived a method to approximate a regularized iterative method
inside $\mathcal{L}$. Starting with $\boldsymbol{y}^0 = \boldsymbol{0}$, we use
the following iterations:
\begin{equation}
  \begin{aligned}
  \boldsymbol{x}_s^{k+1} &= \mathit{FBP}_{\mathcal{L}}(\boldsymbol{p},\boldsymbol{u}_{k+1})  = \boldsymbol{W}_{\mathcal{L}}^T \boldsymbol{C}_{\boldsymbol{u}_{k+1}} \boldsymbol{p} \\
  \boldsymbol{y}^{k+1} &= \boldsymbol{y}^{k} - \alpha \boldsymbol{W}_{\mathcal{L}}^T \boldsymbol{W}_{\mathcal{L}} \boldsymbol{y}^{k} + \boldsymbol{\mathit{d}}^{k+1}\\
  \boldsymbol{x}^{k+1} &= \boldsymbol{x}_s^{k+1} + \boldsymbol{y}^{k+1}
\end{aligned}
\label{eq:localmethod}
\end{equation}
Note that every projection operation in \cref{eq:localmethod} is local, and can
therefore be computed efficiently. The needed filters $\boldsymbol{u}_{k}$ for
all iterations can be precomputed for a certain projection geometry with a
single run of \cref{alg:filtcalc} by returning a filter for each iteration. The
method is based on three approximations to a standard regularized iterative method:
\begin{enumerate}
  \item Iterations of SIRT are approximated by FBP with specific filters.
  \item The prior knowledge is only applied inside $\mathcal{L}$.
  \item The effect of the local prior on pixels outside $\mathcal{L}$ is ignored.
\end{enumerate}
Results from \cref{sec:res} will show that despite these approximations,
reconstructions computed by our method are of significantly higher quality than
either local FBP or global SIRT reconstructions, and visually similar to
global regularized iterative reconstructions. The method is summarized in
\cref{alg:method}.

\begin{algorithm}[t]
\begin{algorithmic}
  \REQUIRE $\boldsymbol{p} \in \mathbb{R}^{N_d N_{\theta}}$, $\boldsymbol{W} \in \mathbb{R}^{N_d N_{\theta} \times N^2}$, $n \in \mathbb{Z}^+$, $\alpha \in \mathbb{R}$
  \STATE $\boldsymbol{y}^0 \leftarrow \boldsymbol{0}$
  \FOR{$k=1$ \TO $n$}
  \STATE $\boldsymbol{x}_s^k \leftarrow \mathit{FBP}_{\mathcal{L}}(\boldsymbol{p},\boldsymbol{u}_k)$
  \STATE $\boldsymbol{y}^k \leftarrow \boldsymbol{y}^{k-1} - \alpha \boldsymbol{W}_{\mathcal{L}}^T \boldsymbol{W}_{\mathcal{L}} \boldsymbol{y}^{k-1} + \boldsymbol{\mathit{d}}^k$
  \ENDFOR
  \RETURN $\boldsymbol{x}_s^n + \boldsymbol{y}^n$
\end{algorithmic}
\caption{Compute a local approximation to a regularized iterative method}
\label{alg:method}
\end{algorithm}

The term $\boldsymbol{\mathit{d}}$ in \cref{alg:method} is the term in which
the prior knowledge is exploited, and depends on which regularized iterative method is
used. Often, in actual implementations, a different formulation can be used that
is more natural to that specific regularized iterative method than the one shown in
\cref{alg:method}. As an example, \cref{alg:methodfista} shows an implementation
of the method when using FISTA to minimize the $\ell_1$ norm of the gradient of
the reconstructed image. Here, we use similar notation to \cite{Beck2009a}, and
$\mathit{FGP}(\boldsymbol{x}, n_{\mathit{FGP}})$ refers to the FGP method of
\cite{Beck2009a} with $n_{\mathit{FGP}}$ iterations, applied to the image
$\boldsymbol{x}$.

\begin{algorithm}[t]
{\fontsize{10}{14}\selectfont
\begin{algorithmic}
  \REQUIRE $\boldsymbol{p} \in \mathbb{R}^{N_d N_{\theta}}$, $\boldsymbol{W} \in \mathbb{R}^{N_d N_{\theta} \times N^2}$, $n \in \mathbb{Z}^+$, $n_{\mathit{FGP}} \in \mathbb{Z}^+$, $\alpha \in \mathbb{R}$
  \STATE $t^0 \leftarrow 1$
	  \STATE $\boldsymbol{x}_{\mathcal{L}}^0 \leftarrow \boldsymbol{0}$
  \STATE $\boldsymbol{x}^0 \leftarrow \boldsymbol{0}$
  \FOR{$k=1$ \TO $n$}
    \STATE $\boldsymbol{x}_s \leftarrow \mathit{FBP}_{\mathcal{L}}(\boldsymbol{p},\boldsymbol{u}_k)$
    \STATE $\boldsymbol{q} \leftarrow \boldsymbol{x}_{\mathcal{L}}^{k-1} - \alpha \boldsymbol{W}_{\mathcal{L}}^T \boldsymbol{W}_{\mathcal{L}} \boldsymbol{x}_{\mathcal{L}}^{k-1}$
    \STATE $\boldsymbol{x}^k \leftarrow \mathit{FGP}(\boldsymbol{x}_s + \boldsymbol{q}, n_{\mathit{FGP}})$
    \STATE $t^k \leftarrow (1+\sqrt{1 + 4 t^{k-1}})/2$
	    \STATE $\boldsymbol{r} \leftarrow \boldsymbol{x} + (t^{k-1}-1) \boldsymbol{x}^{k}/(t^k \boldsymbol{x}^{k-1})$
    \STATE $\boldsymbol{x}_{\mathcal{L}}^k \leftarrow \boldsymbol{r} - \boldsymbol{x}_{s}^{k}$
  \ENDFOR
  \RETURN $\boldsymbol{x}^n$
\end{algorithmic}
}
\caption{Compute a local approximation to FISTA minimizing $\lVert \nabla \boldsymbol{x} \rVert_1$}
\label{alg:methodfista}
\end{algorithm}

\subsection{Implementation details}
\label{sec:impldet}
In this section, we will discuss a few details on implementing the
proposed method. Specifically, we will discuss how to prevent certain
reconstruction artifacts from appearing and how to improve the computation
time of the method in repeated applications.

Using some forms of prior knowledge, artifacts can appear in the reconstructed
image near the edges of the reconstruction grid. For example, the gradient in
a total variation constraint is often defined differently for pixels on the edge
of the reconstruction grid compared to pixels in the interior, which can affect
the reconstruction near the edges. For global regularized iterative methods, the interesting
features of the reconstructed object are usually situated relatively far from the edge,
in which case the artifacts near edges can simply be ignored. In the proposed
local method, however, interesting features may be located near or on
the edge of the chosen local part. A simple but effective way of reducing the
effect of edge artifacts in these cases is to increase the size of the local
part slightly, and crop the resulting reconstruction to the chosen local part.
In the rest of this paper, we increase the size of the local part by padding it
with $\frac{1}{8}$ of the height/width of the local part on each side.

The reconstruction quality of the filter-based approximation of the SIRT
method given in \cref{sec:appsirt} depends on the discrete implementations of the projection operators, as
explained in~\cite{Pelt2015}. Specifically, the method is based on
approximating the combined $\boldsymbol{W}^T\boldsymbol{W}$ operator by a shift-invariant
convolution operation. The discrete projection operations can be implemented in different
ways~\cite{Xu2006}, and the accuracy of the approximation depends on the
chosen implementation. In practice, most artifacts resulting from the
errors in the approximation
are found in the low frequencies of the reconstructed image, similar to the
artifacts that can occur when discretizing the Ram-Lak filter of the FBP method~\cite[Fig. 3.13]{Kak2001}.
By using implementations of the projection operators that minimize the approximation error
that is made, reconstruction artifacts can be limited, and are typically invisible to a human
observer.
In this paper, we use an additional preprocessing step to further reduce these artifacts. Before each
reconstruction with the local approximation method,
we subtract from the projection data the forward projection of a disc,
centered on the rotation axis, with a diameter $N$ and a constant gray value. The gray value is chosen such that the
$\ell_2$-norm of the zero-frequency components of all projections are minimized after subtraction.
By reducing the low-frequency components of the projection data with this procedure,
the artifacts resulting from the approximation error are reduced as well.
After reconstruction,
the same disc is added back to the reconstructed image. In practice, this procedure ensures that artifacts resulting
from errors made in approximating SIRT by filtered backprojection are minimal.

As explained in \cref{sec:prior}, all projection operations of the proposed method
can be computed locally, and are therefore efficient to compute. When the local
part is much smaller than the number of detector pixels ($N_{\mathcal{L}}\ll
N_d$), however, the convolution operation in $\mathit{FBP}_{\mathcal{L}}$,
which scales with $N_d$ instead of $N_{\mathcal{L}}$, can become
a significant part of the total computation time. In many cases, however, one
will perform repeated applications of the local method, for example
when finding optimal parameters for the applied prior knowledge term, or when
reconstructing multiple local parts at different locations. In these cases, the
convolution of the projection data with the different filters $\boldsymbol{u}_k$
for each iteration can be precomputed
once and reused for the different local reconstructions, improving
reconstruction time significantly.

\section{Experiments}
\label{sec:exp}
To investigate the properties of the proposed method, we implemented it in
Python, version 3.4.3, using the ASTRA toolbox~\cite{Aarle2015} to perform all
tomographic projection operations, which enables the use of optimized GPU-based
computations~\cite{Palenstijn2011}. All experiments were performed on a machine
running Fedora Linux 21, with an Intel Xeon E5-2623 processor, 13 GB of memory,
and a NVIDIA GeForce GTX TITAN Z GPU using CUDA version 7.\@0.

We present results for three different forms of prior knowledge about the
reconstructed object: one domain constraint and two penalty constraints.
For the domain constraint we use box constraints on the pixel values
by specifying $D=\{\boldsymbol{x} \in \mathbb{R}^{N^2}; l \leq x_i \leq r, i=1,\dots,N^2\}$
in the objective function of~\cref{eq:priorsys}.
For the penalty constraints, we use  $\ell_1$ minimization
of the reconstruction in a Haar wavelet basis, i.e. specifying $g(\boldsymbol{x}) = \lVert \boldsymbol{B}\boldsymbol{x}\rVert_1$,
and $\ell_1$ minimization of the
gradient of the reconstructed image (TV minimization), i.e. specifying $g(\boldsymbol{x}) = \lVert \nabla \boldsymbol{x}\rVert_1$. We use~\cref{eq:sirtbox} to
find solutions in the case of box constraints on the pixel values, and the FISTA
method in the case of both $\ell_1$ penalty functions. In all cases, we
compare the locally approximated reconstructions with global reconstructions of the
full object exploiting the same prior knowledge on the full volume, and with the popular analytical
FBP method and algebraic SIRT method, which do not explicitly exploit any prior knowledge.

\begin{figure}
	\begin{center}
		\subfloat[][]{\frame{\includegraphics[width=\figwidth]{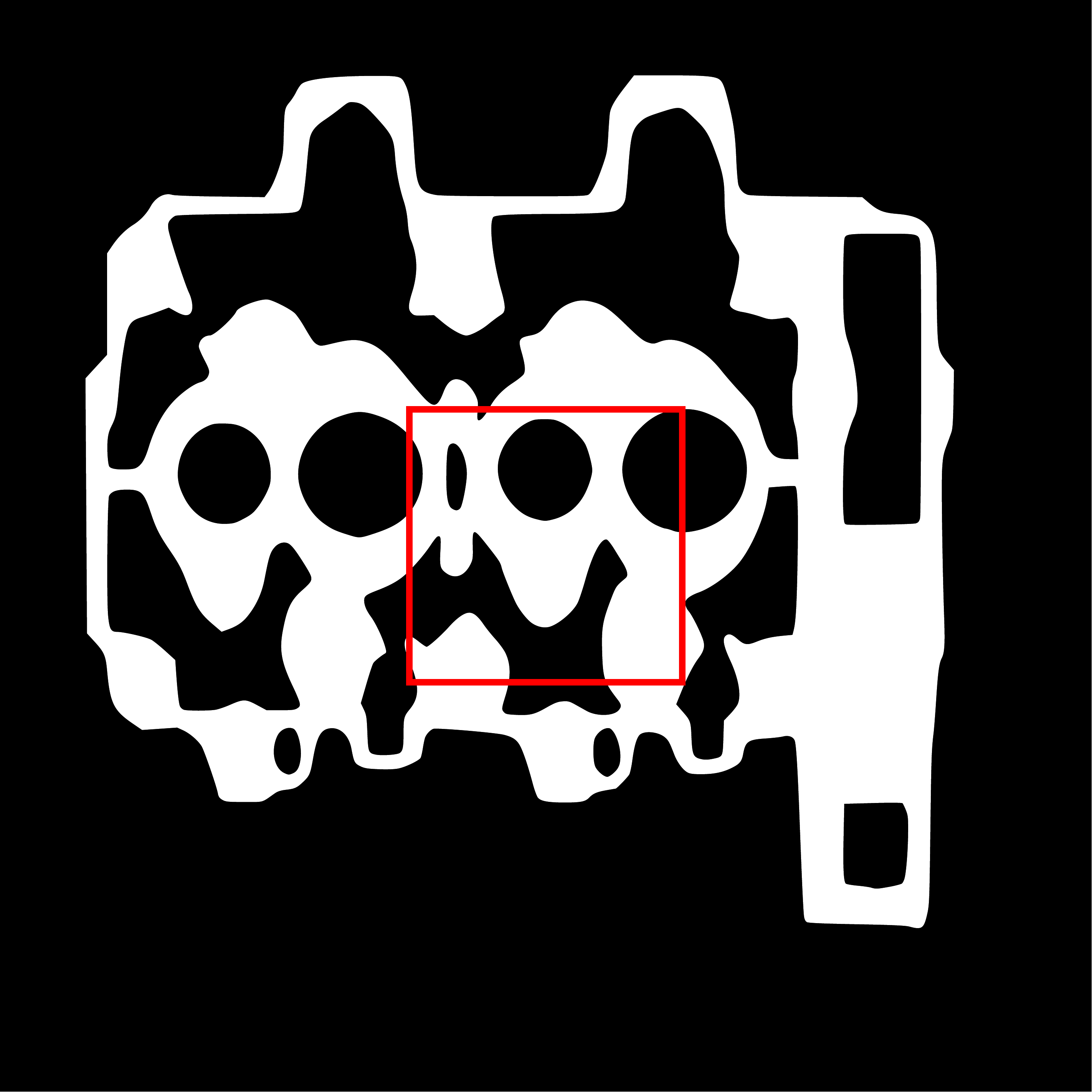}}}
		\hfill
\subfloat[][]{\frame{\includegraphics[width=\figwidth]{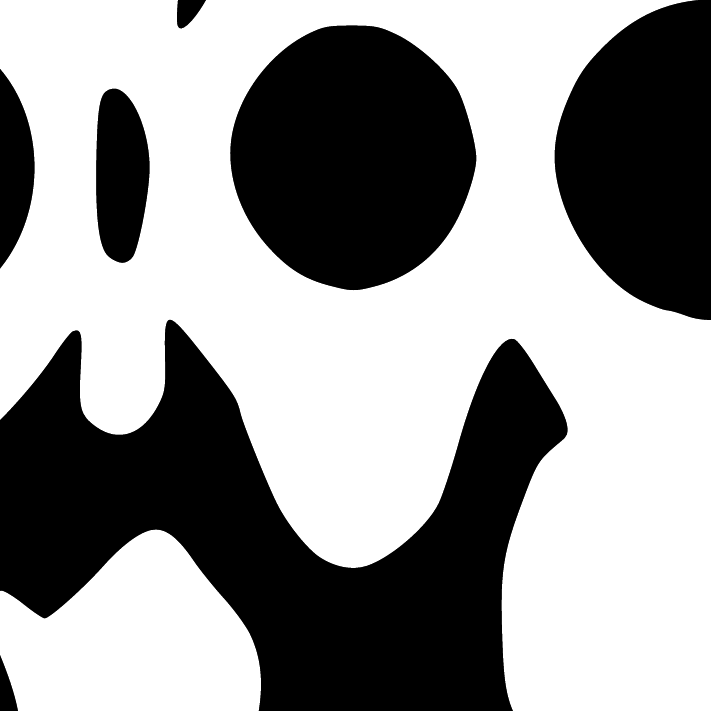}} \label{fig:motorfull}}
		\hfill
\subfloat[][]{\frame{\includegraphics[width=\figwidth]{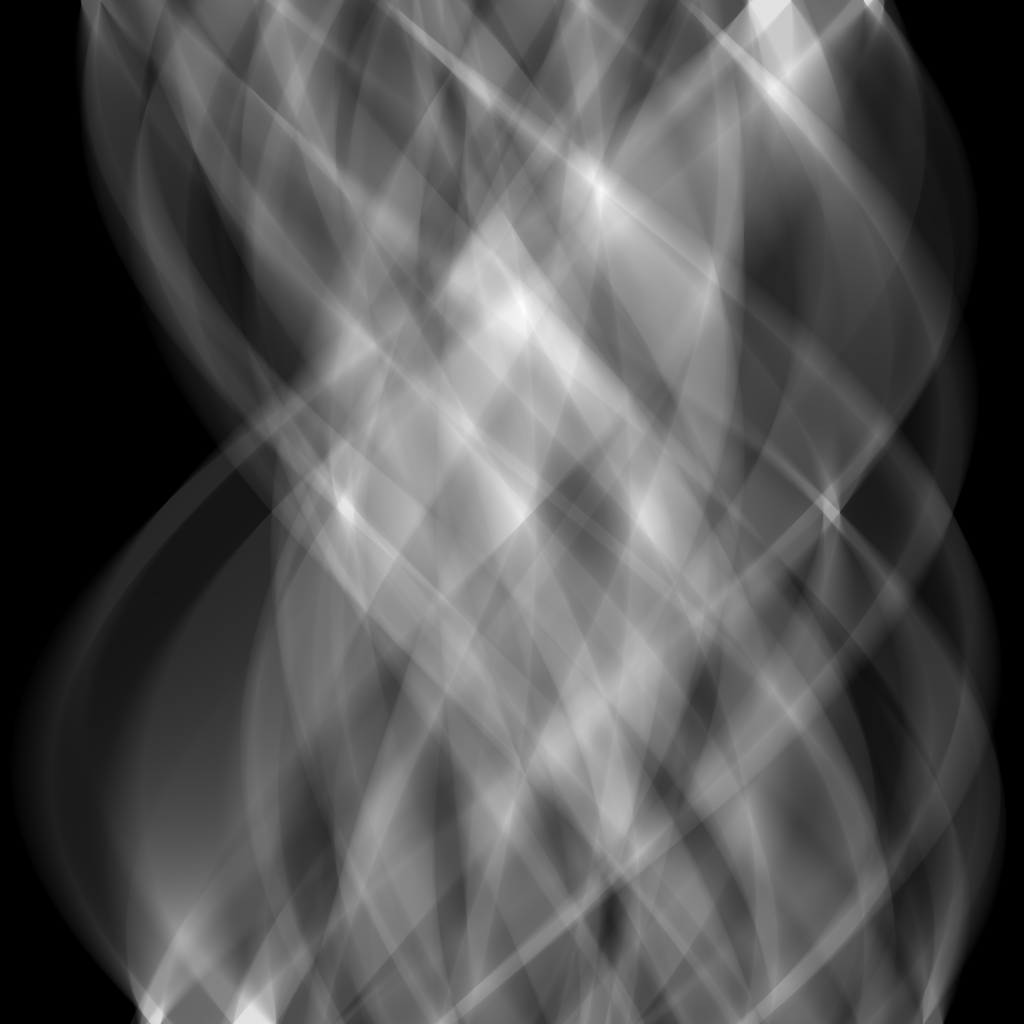}}}
	\end{center}
	\caption{The phantom used for most experiments in this paper. In (a), the
		entire phantom is shown with a red square indicating the local
	part (b) that is used in most experiments. In (c), the sinogram of the
phantom is shown for 1024 detector pixels and 1024 projections equally
distributed in $[0,\pi]$.}
\label{fig:phantom}
\end{figure}

The phantom that is used in most experiments in this paper is shown in \cref{fig:phantom}. This phantom
was chosen because it is suitable for all three forms of prior knowledge that we exploit. It consists
of two materials: a background with a value of zero and a foreground with a value of one. Therefore,
box constraints can be effectively exploited by setting $l=0$ and $r=1$.
Since the phantom has a sparse boundary,
TV minimization and a Haar wavelet basis can also be used to improve
reconstruction quality.
In addition to the phantom shown in \cref{fig:phantom}, we also present some results
for the Shepp-Logan head phantom, shown in~\cref{fig:slphan},
which has a relatively sparse boundary as well.

For each reconstruction, we report the mean squared error (MSE) of the
reconstructions inside the region of interest, compared to a known ground truth
image. We also report the structural similarity index (SSIM) \cite{Wang2004} of
the reconstructions inside the region of interest compared to the ground truth, which is a metric that is
designed to be closer to the human visual system than the mean squared error. For
methods where a parameter needs to be chosen, i.e. $\lambda$ in
\cref{eq:priorsys}, we perform two reconstructions each time: one with the value that
minimizes the $\mathit{MSE}$ and one with the value that maximizes the
$\mathit{SSIM}$. In each case, we find the optimal parameter value using the
Nelder-Mead method~\cite{Nelder1965}. Note that the optimal parameter value can
depend on the dimensions of the reconstruction grid, and therefore, the optimal values
can be different for the global regularized iterative reconstructions compared to the locally
approximated reconstructions. For all iterative methods, we use 200 iterations to compute
each reconstruction, and we use 100 FGP iterations in the FISTA
method for TV minimization~\cite{Beck2009a}.

In most experiments, we use a $4096\times4096$ pixel image of the phantom, and generate
projection data for $4096$ detector pixels. Afterwards, the projection data is resampled
to $1024$ detector pixels, and reconstructions are computed on a $1024\times1024$ pixel
grid, or a local part of that grid. These reconstructions are compared to the original
$4096\times4096$ pixel phantom, resampled to a $1024\times1024$ pixel grid. In most cases,
additional Poisson noise is applied to the projection data to simulate experimental
conditions. The amount of applied Poisson noise is indicated by a variable $I_0$, with
lower values corresponding to higher amounts of applied noise. Specifically, the noise
is applied by first transforming the simulated projections to virtual photon counts,
in which the largest photon count out of all detector pixels is set to $I_0$. For each
detector pixel, a new photon count is sampled from a Poisson distribution with the original
photon count as the expected value. Finally, the resulting noisy photon counts are transformed
back to noisy line integrals of the phantom.

\section{Results}
\label{sec:res}
In this section, we present the results of the experiments that we performed to
investigate the properties of the proposed local approximation method, and discuss
these results.
\subsection{Local SIRT approximation}
\label{sec:locsirtres}
\begin{figure}
\begin{center}
\subfloat[][]{\frame{\includegraphics[width=\figwidth]{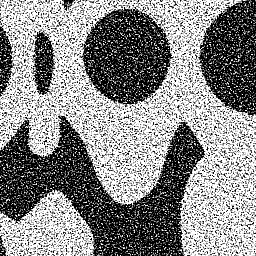}}}
\hfill
\subfloat[][]{\frame{\includegraphics[width=\figwidth]{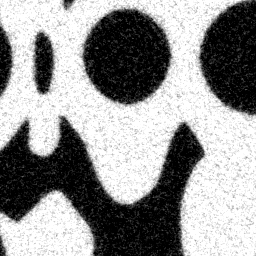}}}
\hfill
\subfloat[][]{\frame{\includegraphics[width=\figwidth]{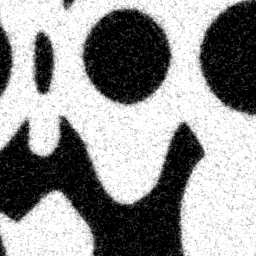}}}
\end{center}
\caption{Reconstructions of a $256\times256$ pixel local part of the motor phantom
(\cref{fig:motorfull}), using projection data of $1024$ detector pixels with
$N_\theta=512$ projection angles, equally distributed in the interval $[0,\pi]$,
and with Poisson noise applied. In (a) the local FBP reconstruction is shown, in
(b) the global SIRT reconstruction cropped to the local part, and in (c) the locally
approximated SIRT reconstruction.
}
\label{fig:locsirt}
\end{figure}

In \cref{fig:locsirt}, reconstructions are shown for the local part of the phantom,
computed by standard FBP, standard SIRT, and the local approximation of SIRT
(\cref{eq:locsirtapp}). Note that the global SIRT reconstruction and its local
approximation are visually very similar. The difference between the computation
times is significant, however: the local reconstructions take 28 milliseconds to
compute each, while the global SIRT reconstruction takes 2.6 seconds. The $\mathit{MSE}$
of the FBP, SIRT, and local approximation are 0.245, 0.016, and 0.016, respectively, and
the $\mathit{SSIM}$ values are 0.07, 0.25, and 0.27.

\subsection{Local regularized iterative approximation}

\begin{figure}
\begin{center}
\includegraphics[width=\linewidth]{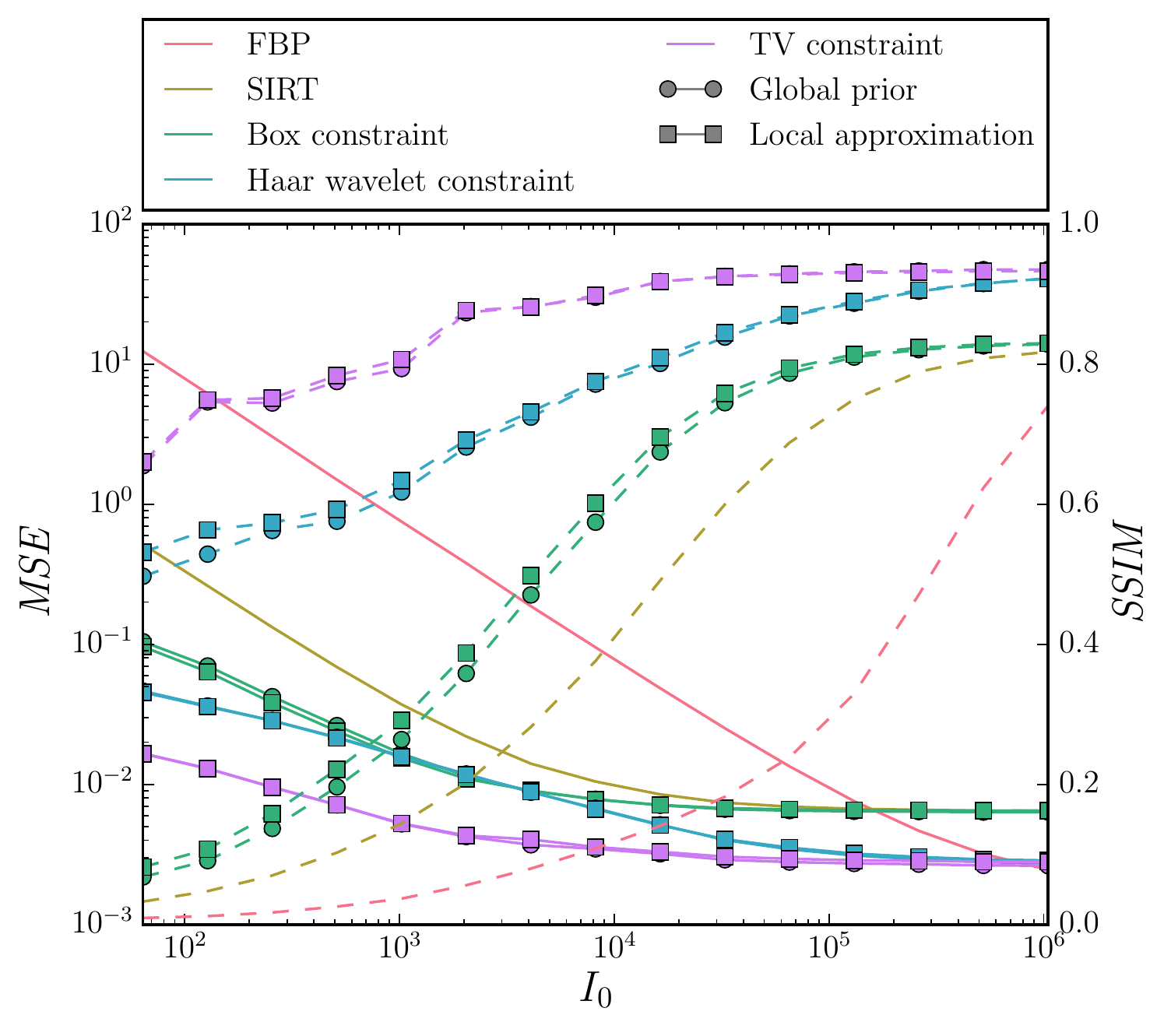}
\end{center}
\caption{Mean squared error ($\mathit{MSE}$, solid lines) and structural similarity index
	($\mathit{SSIM}$, dashed lines) of reconstructions of a region ($256\times256$
pixels) of the motor phantom (\cref{fig:motorfull}) for various amounts of applied
Poisson noise $I_{0}$ and types of prior knowledge.
The reconstructions are computed using projection data of 1024 detector pixels and 512 projections equally
distributed in the interval $[0,\pi]$.}
\label{fig:noises}
\end{figure}

\begin{figure}
\begin{center}
\includegraphics[width=\linewidth]{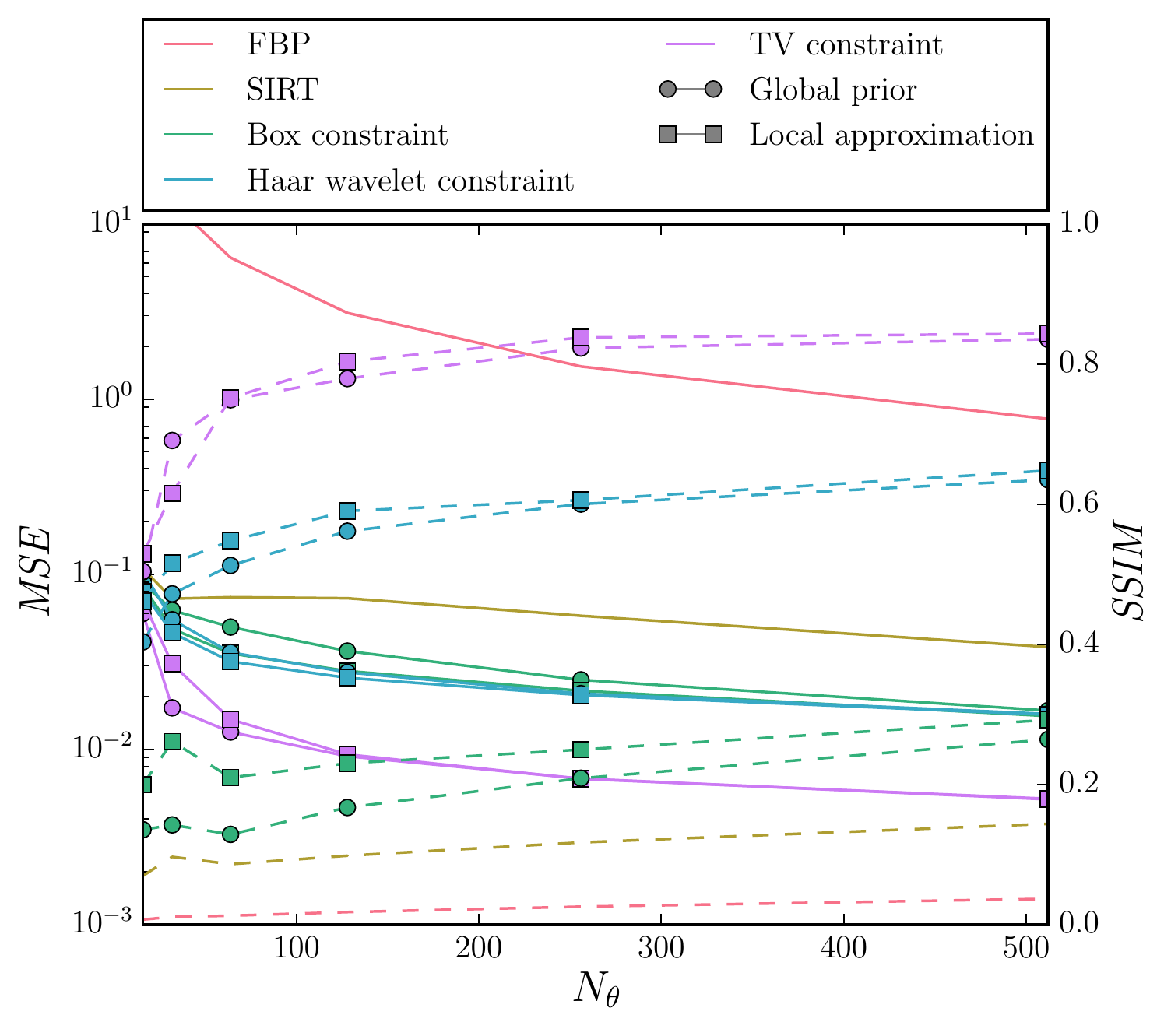}
\end{center}
\caption{Mean squared error ($\mathit{MSE}$, solid lines) and structural similarity index
	($\mathit{SSIM}$, dashed lines) of reconstructions of a region ($256\times256$
pixels) of the motor phantom (\cref{fig:motorfull}) for various numbers of projection angles $N_{\theta}$
(equally distributed in the interval $[0,\pi]$)
and types of prior knowledge.
The reconstructions are computed using projection data of 1024 detector pixels, with
applied Poisson noise.}
\label{fig:angles}
\end{figure}

In \cref{fig:noises}, the mean squared error and structural similarity index
are shown as a function of the amount of applied
Poisson noise $I_{0}$, for standard FBP, standard SIRT, and global and locally
approximated reconstructions using various types of prior knowledge. The results
show that by exploiting prior knowledge, reconstruction quality can be significantly
improved compared to standard FBP and SIRT reconstructions. For this phantom,
exploiting total variation minimization yields reconstructions with the
lowest $\mathit{MSE}$ and highest $\mathit{SSIM}$ values. The results also show that
for all tested types of prior knowledge, the quality metrics of the locally approximated
reconstructions are very close to those of the global regularized iterative reconstructions. For unknown reasons,
the quality metrics of the local approximations are
slightly better than the global regularized iterative reconstructions.
Similar results can be seen in \cref{fig:angles}, where the quality metrics are shown as a function of the number of projections
angles.

\begin{figure}
	\begin{center}
\includegraphics[width=\linewidth]{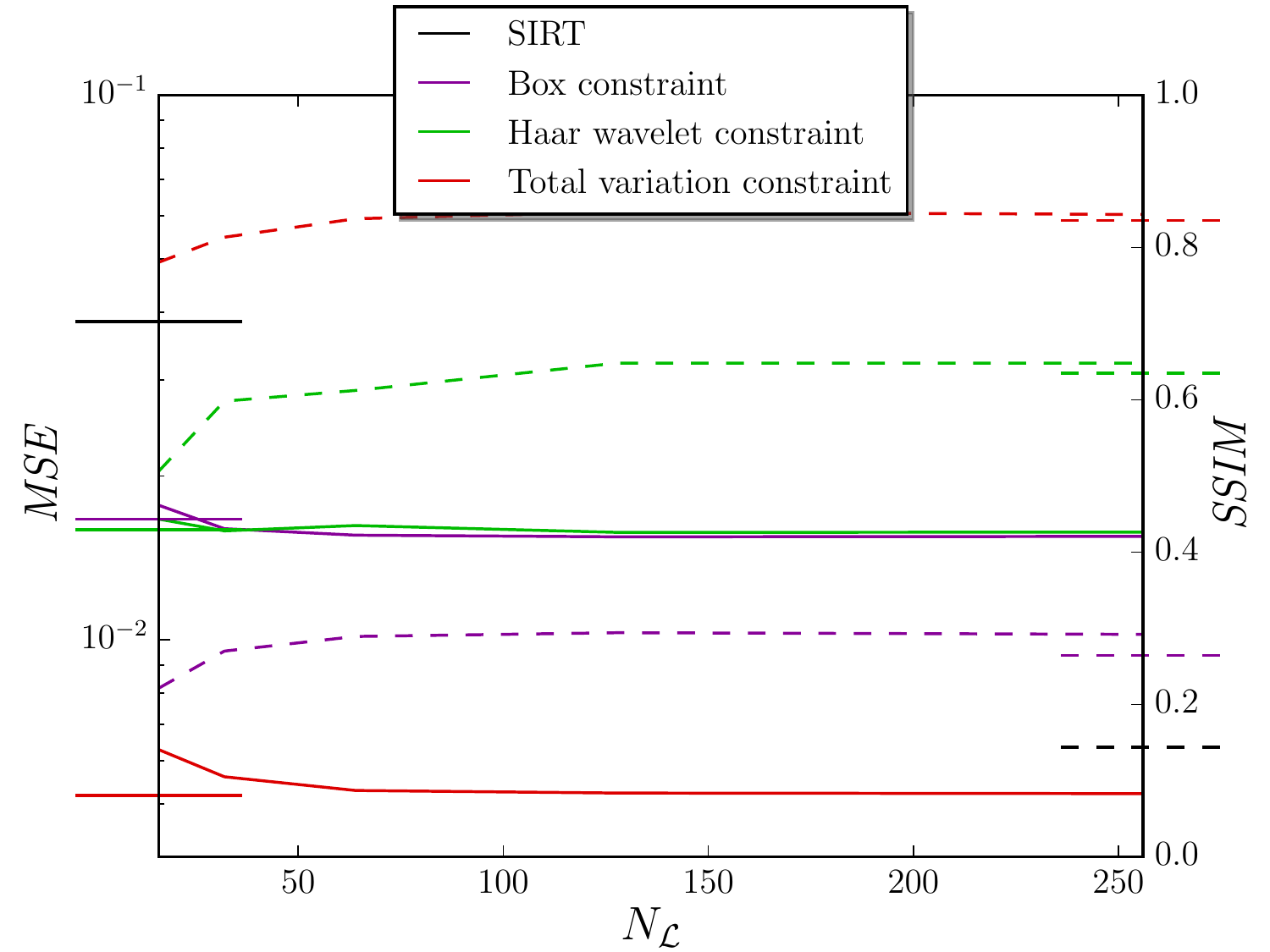}
	\end{center}
\caption{Mean squared error ($\mathit{MSE}$, solid lines) and structural similarity index
	($\mathit{SSIM}$, dashed lines) of reconstructions of a region ($256\times256$
pixels) of the motor phantom (\cref{fig:motorfull}) for various sizes of the local
part $N_{\mathcal{L}}$ and types of prior knowledge.
The reconstructions are computed using projection data of 1024 detector pixels and 512 projections equally
distributed in the interval $[0,\pi]$, with
applied Poisson noise. For $N_{\mathcal{L}}<256$, multiple local reconstructions are tiled to
create a reconstruction of $256\times256$ pixels, to enable comparison between different local sizes.
The partial horizontal lines on each axis indicate the $\mathit{MSE}$ and
$\mathit{SSIM}$ of global SIRT and global regularized iterative reconstructions, cropped to the
same $256\times256$ pixels.
}
\label{fig:sizes}
\end{figure}

The mean squared error and structural similarity index are shown as a function of the size
of the local part $\mathcal{L}$ in \cref{fig:sizes}. For all three prior knowledge types,
the reconstruction quality of the local approximations is only significantly lower
compared to the global regularized iterative methods when the local size is
$N_{\mathcal{L}} = 32$ or smaller, at which point the number of pixels of the local part
is less than 0.1\% of the number of pixels in the global reconstruction grid.
For larger local sizes, the reconstruction quality is almost independent
of the local size.
These results suggest that, even for reasonably small local parts, the approximations that are made by
the proposed local method do not influence the reconstruction quality significantly.

\begin{figure}
\begin{center}
	\subfloat[][]{\frame{\includegraphics[width=\figwidth]{motor.pdf}}
	\label{fig:locphan}}
\hfill
\subfloat[][]{\frame{\includegraphics[width=\figwidth]{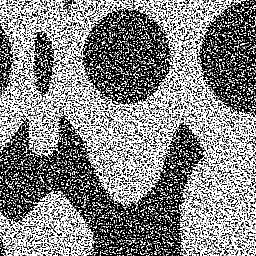}}}
\hfill
\subfloat[][]{\frame{\includegraphics[width=\figwidth]{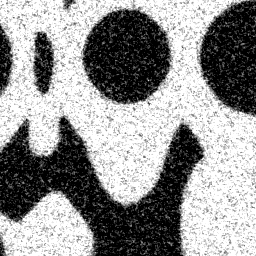}}}
\\
\subfloat[][]{\frame{\includegraphics[width=\figwidth]{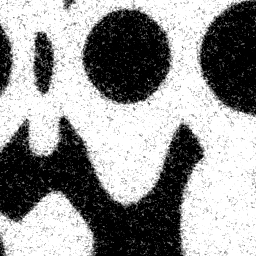}}}
\hfill
\subfloat[][]{\frame{\includegraphics[width=\figwidth]{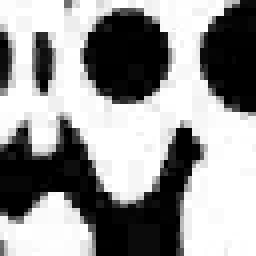}}}
\hfill
\subfloat[][]{\frame{\includegraphics[width=\figwidth]{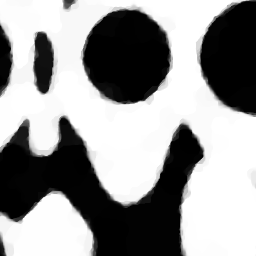}}}
\\
\subfloat[][]{\frame{\includegraphics[width=\figwidth]{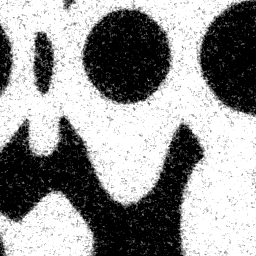}}}
\hfill
\subfloat[][]{\frame{\includegraphics[width=\figwidth]{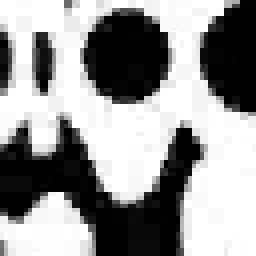}}}
\hfill
\subfloat[][]{\frame{\includegraphics[width=\figwidth]{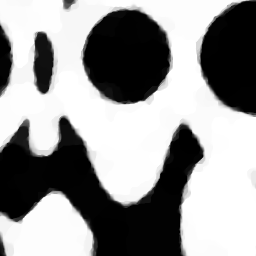}}}
\end{center}
\caption{Reconstructions of a local part ($256\times256$ pixels) of the motor
phantom (a) from projection data of 1024 detector pixels and 512 projections
equally distributed in the interval $[0,\pi]$, with Poisson noise applied, using various reconstruction
methods: (b) local FBP, (c) global SIRT cropped to local part,
(d)-(f) global regularized iterative method cropped to local part, with (d) box constraint,
(e) Haar wavelet constraint, and (f) TV constraint, and (g)-(i) the proposed
local method, with (g) box constraint, (h) Haar wavelet constraint, and
(i) TV constraint.
The local reconstructions are shown with a gray-level window in which black corresponds to the minimum
value and white to the maximum value of the phantom inside the local part.
}
\label{fig:recs}
\end{figure}
\begin{figure}
\begin{center}
    \subfloat[][]{\frame{\includegraphics[width=\figwidth]{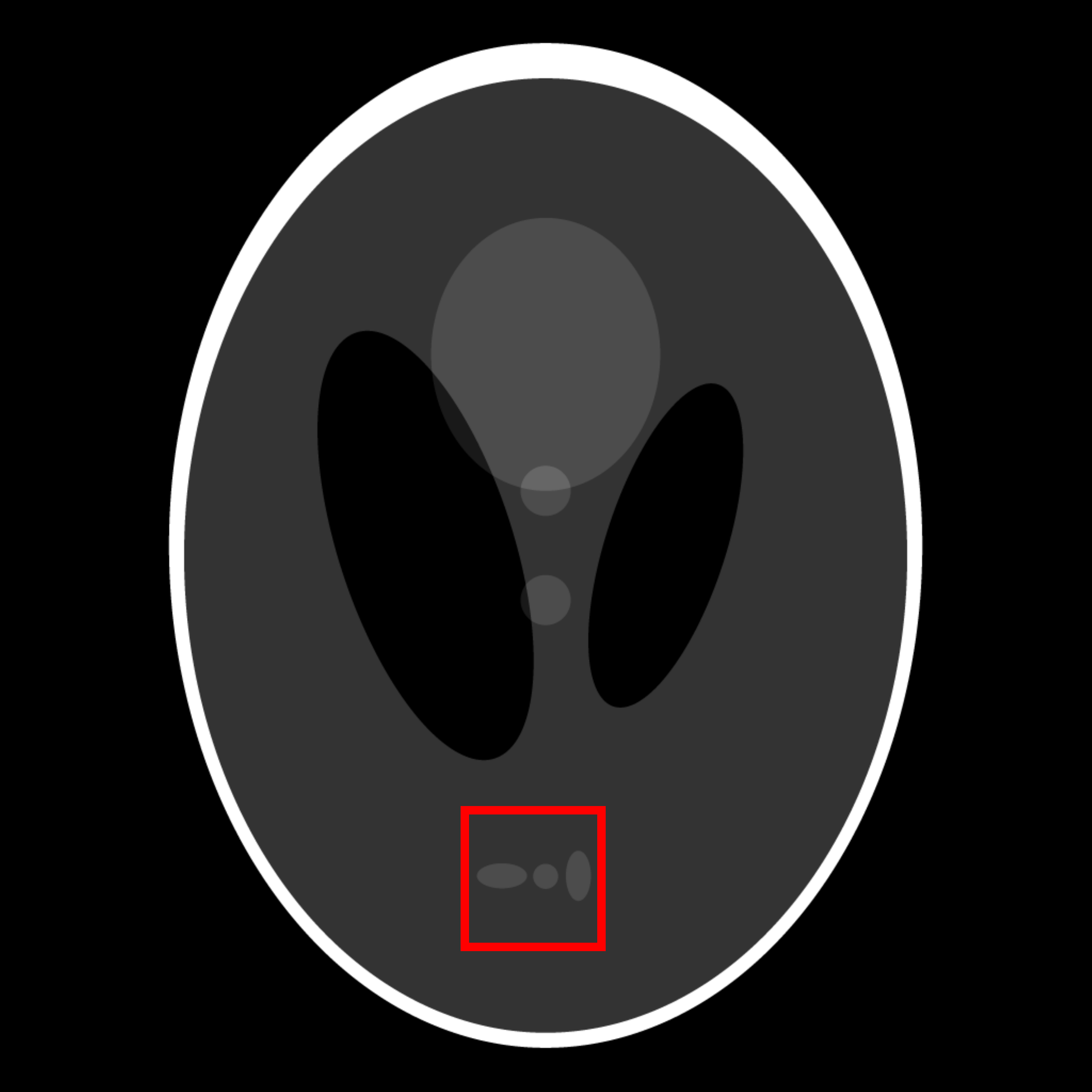}}\label{fig:slphan}}
\hfill
\subfloat[][]{\frame{\includegraphics[width=\figwidth]{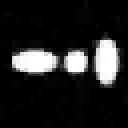}}}
\hfill
\subfloat[][]{\frame{\includegraphics[width=\figwidth]{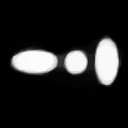}}}
\\
\subfloat[][]{\frame{\includegraphics[width=\figwidth]{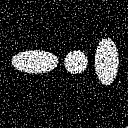}}}
\hfill
\subfloat[][]{\frame{\includegraphics[width=\figwidth]{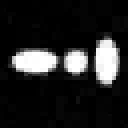}}}
\hfill
\subfloat[][]{\frame{\includegraphics[width=\figwidth]{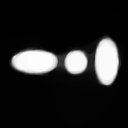}}}
\end{center}
\caption{Reconstructions of a local part ($128\times128$ pixels) of the Shepp-Logan head phantom,
indicated by the red square in (a). The reconstructions are computed from projection data of 1024 detector pixels and 512 projections
equally distributed in the interval $[0,\pi]$, with Poisson noise applied, using various reconstruction
methods:
(b)-(c) global regularized iterative method cropped to local part, with
(b) Haar wavelet constraint and (c) TV constraint, (d) local FBP, and (e)-(f) the proposed
local method, with (e) Haar wavelet constraint and
(f) TV constraint.
The local reconstructions are shown with a gray-level window in which black corresponds to the minimum
value and white to the maximum value of the phantom inside the local part.
}
\label{fig:recssl}
\end{figure}

Reconstructed images of a local part with $256\times 256$ pixels are shown in \cref{fig:recs}, for
projection data of 1024 detector pixels and 512 equiangular projections with Poisson noise applied.
The images show that the local approximations are visually almost identical to the global regularized
iterative reconstructions for all three prior knowledge types. The results also show how the different
prior knowledge types can help improve certain image characteristics compared to standard FBP and SIRT
reconstructions. In \cref{fig:recssl}, reconstructed images are shown for a smaller local part ($128 \times
128$ pixels) of the Shepp-Logan head phantom. Similar to the previous results, the local approximations are
visually almost identical to the global regularized iterative reconstructions.

\subsection{Computation time}
\begin{figure}
	\begin{center}
\includegraphics[width=\linewidth]{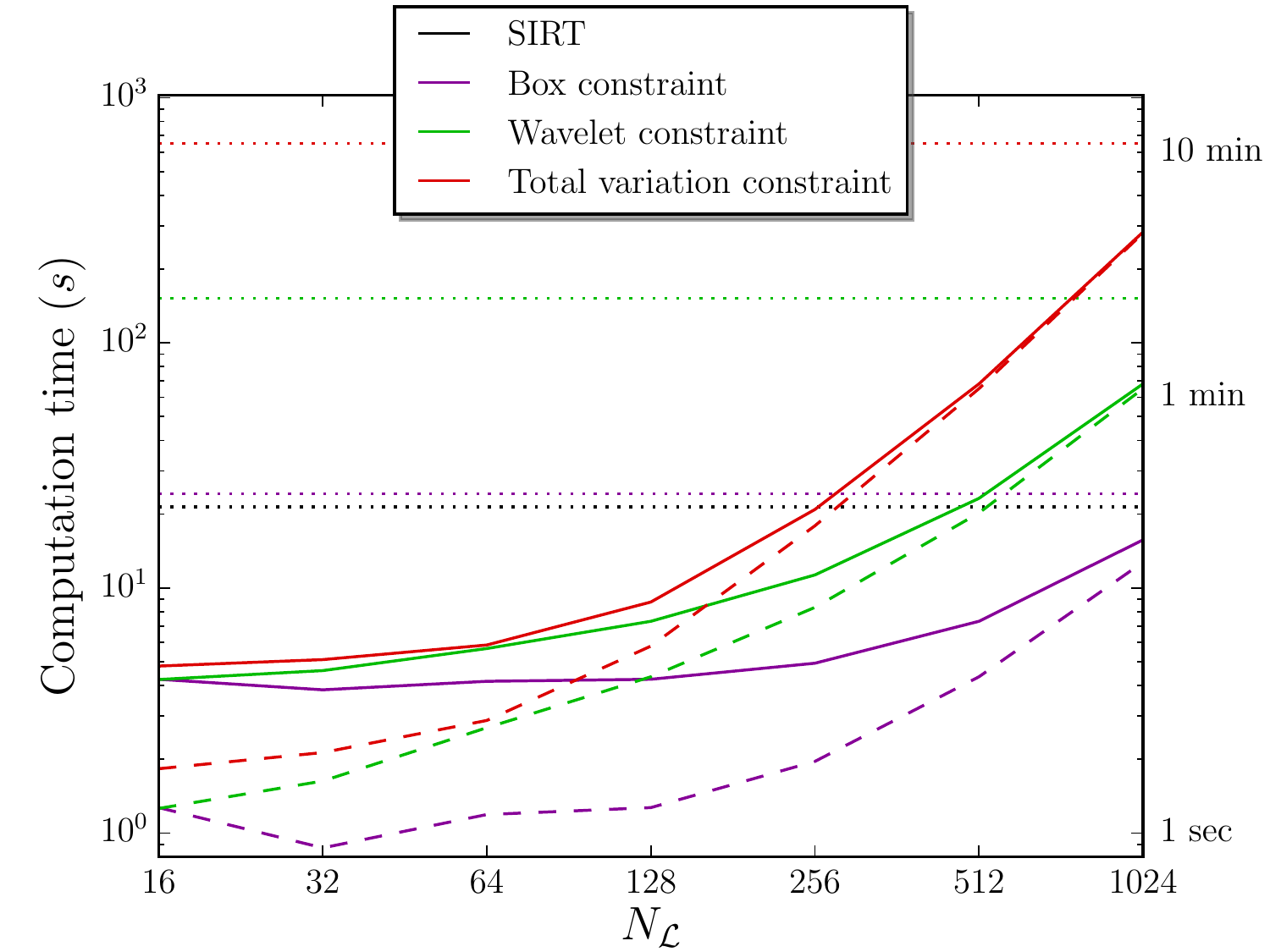}
	\end{center}
\caption{Reconstruction time of the global regularized iterative methods (dotted) and the
	proposed local method for various sizes of the local part
$N_{\mathcal{L}}$ and constraint types, using data of 2048 detector pixels
and 512 projections. Solid lines show the reconstruction time for a single
application of the local method, and dashed lines show the reconstruction
time for subsequent applications, where the convolution results of an earlier
reconstruction can be reused (see \cref{sec:impldet}).}
\label{fig:time}
\end{figure}

The computation time of the proposed local reconstruction method is shown in \cref{fig:time} as a function
of the size of the local part $\mathcal{L}$. Also shown is the computation time of the standard global
regularized iterative method. For the local method, computation times are shown both for the first application, as
well as for subsequent applications, in which the convolution results of the first application can be reused to
decrease the needed computation time (see \cref{sec:impldet}). For all types of prior knowledge, the local method
requires significantly less computation time than the global regularized iterative methods.

If one is only
interested in a local part of the object, the local method can be used to compute advanced regularized reconstructions
in a few seconds instead of the several minutes it costs to compute the global reconstruction. In cases where the
same regularized iterative method is computed multiple times for the same projection data, for example when estimating
the $\lambda$ parameter, the proposed local method requires even less computation time, leading to a significant reduction
of processing time in practice. Finally, since each local reconstruction is independent of the other local reconstructions,
different local parts can be reconstructed in parallel and combined afterwards to compute a larger part of the scanned object
in short time. An example of such a reconstruction is shown in \cref{sec:tiling}.

\subsection{Experimental data}

\begin{figure*}
\begin{center}
	\subfloat[][]{\includegraphics[width=\figwidthfour]{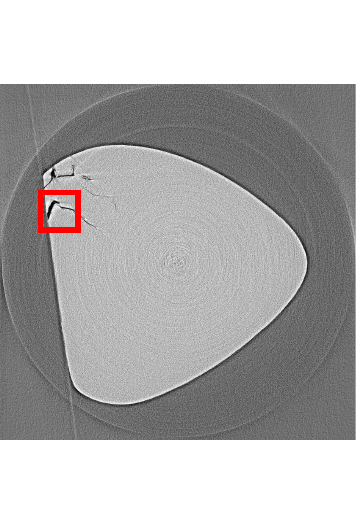}}
	\hfill
	\subfloat[][]{\includegraphics[width=\figwidthfour]{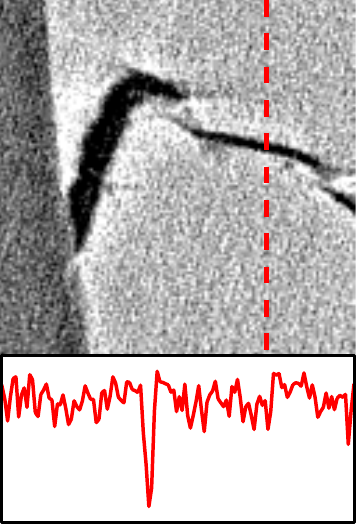}}
	\hfill
	\subfloat[][]{\includegraphics[width=\figwidthfour]{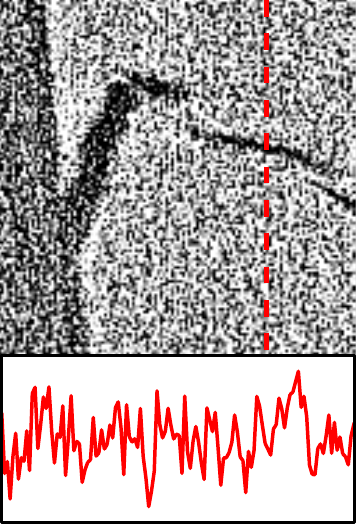}}
	\hfill
	\subfloat[][]{\includegraphics[width=\figwidthfour]{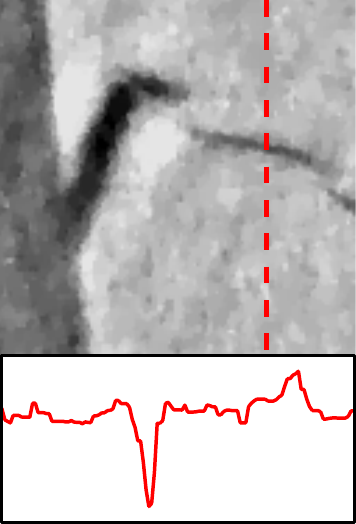}}
\end{center}
\caption{Reconstructions of a local part ($128\times128$ pixels) of experimental data of a small fatigue test sample made
from Ti alloy VST 55531,
acquired with 1200 detector pixels and 1500 projections equally distributed over the interval $[0,\pi]$. In (a) and (b), FBP reconstructions
are shown using all 1500 projections, with the local part indicated by a red square in (a). The local FBP reconstruction
using only 75 equiangular projections is shown in (c). In (d), a reconstruction is shown for the same 75 projections, using the local reconstruction method presented
in this paper with TV-minimization regularization by the FISTA method.
Underneath each local reconstruction, the line profile of the column indicated by the dashed line is shown.
}
\label{fig:exper}
\end{figure*}

In \cref{fig:exper}, reconstructed images are shown for a local part of an experimental
dataset. The experimental data was acquired for a small fatigue test sample made from Ti
alloy VST 55531. The sample was scanned at beamline
ID11 of the European Synchrotron Radiation Facility (ESRF),
with a parallel, monochromatic (52 keV) synchrotron X-ray beam.
The distance between the sample and detector was 40 mm,
and 1500 projections were acquired, equally distributed in the interval $[0,\pi]$.
The projections were acquired on a high resolution detector system,
resulting in projections, after $2 \times 2$ binning, with $1200\times 1200$ pixels and an
effective pixel size of 0.56 microns.

Results are shown in \cref{fig:exper} for a single slice of the reconstructed dataset,
computed using FBP and the proposed local method with a TV minimization constraint. For
FBP, we show results both when using all 1500 projections that were acquired, and when using
only 75 projections, selected by taking every 20\textsuperscript{th} projection of the
full dataset. For the local method, we show results for the same limited dataset of 75 projections.
The results show that the local method can be successfully applied to an experimental dataset to exploit
prior knowledge in the reconstruction. Compared to the FBP reconstruction using 75 projections, the
local method is able to more clearly separate the formed crack from the sample itself, which is
especially visible in the line profiles. Note that in this type of sample, a user would typically
only be interested in the highly localized crack that is forming in the sample, which would make
global regularized iterative methods waste significant amounts of computation time on parts
of the sample that are not interesting. With the proposed local method, on the other hand,
a user would be able to select and reconstruct only those parts of the sample that are interesting.

\subsection{Tiling reconstructions}
\label{sec:tiling}
\begin{figure}
\begin{center}
\subfloat[][]{\includegraphics[width=\figwidth]{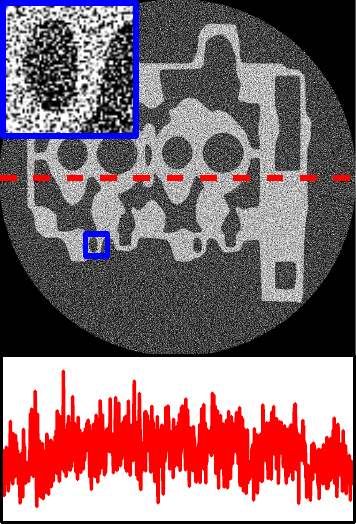}}
\hfill
\subfloat[][]{\includegraphics[width=\figwidth]{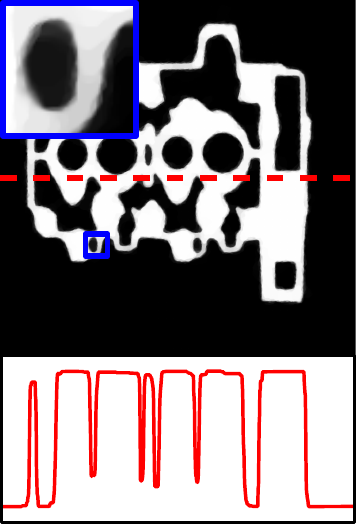}}
\hfill
\subfloat[][]{\includegraphics[width=\figwidth]{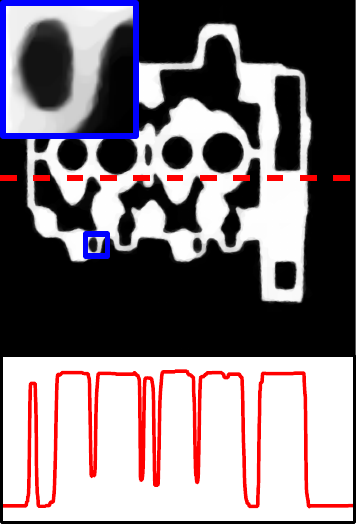}}
\end{center}
\caption{Reconstructions of the motor phantom using projection data of
$1024$ detectors and $N_\theta=512$ projection angles, equally distributed in the interval $[0,\pi]$,
with Poisson noise applied.
The reconstructions are computed with (a) FBP, (b) global TV minimization by the FISTA method, and
(c) local $128\times128$ pixel reconstructions tiled to the complete $1024\times1024$ pixel grid. The local
reconstructions in (c) are computed using the local reconstruction method presented in this paper with TV-minimization regularization by
the FISTA method.
Underneath each reconstruction, the line profile of the row indicated by the dashed line is shown.
A small region, indicated by the blue square, is shown enlarged in the top-left
corner of each reconstruction as well.
}
\label{fig:tile}
\end{figure}
As explained before, one possibility of the proposed local method is to reconstruct different local parts
of the image and combine them afterwards into a single reconstruction. One application of this approach would
be to compute the different local parts in parallel, which can be parallelized efficiently since each local
reconstruction is independent of the others. Another application would be to estimate reconstruction
parameters such as the $\lambda$ term of \cref{eq:priorsys} only in a local part of the reconstruction,
which would significantly reduce the time needed to estimate them. Afterwards,
the complete image can be reconstructed by combining several local reconstructions using these parameters, which
can be computed in parallel as well.

An example of a reconstruction that is computed by tiling several local reconstructions is shown in \cref{fig:tile}.
In this case, we combined 64 local reconstructions of $128\times 128$ pixels each to compute a single $1024 \times 1024$
pixel reconstruction, using TV-minimization as the prior knowledge term. The local reconstructions are tiled by simply
placing them next to each other on the large reconstruction grid, without any overlapping regions.
The results show that
there are no visible artifacts from this tiling procedure.
Furthermore, the tiled reconstruction is visually
almost identical to a reconstruction computed by the global regularized iterative method.
This shows that it is possible to
significantly reduce the computation time of a global regularized iterative reconstruction method by
approximating it with a tiling of local reconstructions computed in parallel.

\subsection{Truncated projection data}
\begin{figure}
\begin{center}
\subfloat[][]{\includegraphics[width=\figwidth]{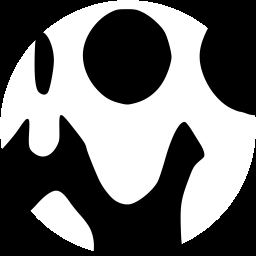}}
\hfill
\subfloat[][]{\includegraphics[width=\figwidth]{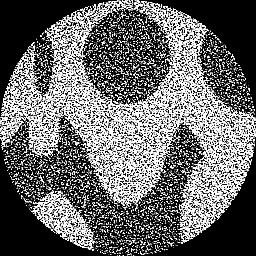}}
\hfill
\subfloat[][]{\includegraphics[width=\figwidth]{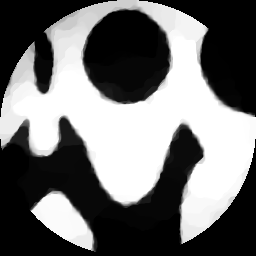}}
\end{center}
\caption{Reconstructions of the motor phantom (a), using projection data of
$1024$ detectors truncated to the central $256$ detector pixels, using
$N_\theta=512$ projection angles, equally distributed in the interval $[0,\pi]$,
and with Poisson noise applied. The reconstructions of FBP (b) and the proposed
method with TV-minimization by the FISTA method (c) are shown for the central disc with a width of
$256$ pixels. Constant padding is used in both reconstructions to reduce truncation artifacts.
}
\label{fig:trunc}
\end{figure}

In some applications of tomography, it is impossible
to acquire projections that include the entire
scanned object.
In these cases, the acquired projection data are
\emph{truncated} at the edge of the detector.
The resulting
reconstruction problem is similar to local reconstruction: again, one is only interested in
a subvolume of the entire scanned object. In this case, however, data for the
object outside the subvolume is missing.
Filtered backprojection is often used to reconstruct truncated data
by simply padding the acquired data in order to reduce the
artifacts caused by the truncation. Since the local method proposed in this paper
uses FBP to approximate the SIRT method, the same padding approach can be used
to apply the method to truncated data.
Reconstructions of truncated phantom data are shown in \cref{fig:trunc}, for
FBP and the proposed local method. The results show that the local method
can be used to exploit prior knowledge in the case of truncated data to improve
reconstruction quality.

\section{Conclusions}
\label{sec:conc}

In this paper, we introduced a method to approximate regularized iterative
tomographic reconstruction methods inside a region of interest. This
method can be used to reduce computation time when one is only interested
in the reconstruction inside the region of interest. The method
is based on approximating the SIRT steps that are part of many
regularized iterative methods by filtered backprojection with specific
pre-calculated filters. The result is a reconstruction method in which
all projection operations involve only the pixels that are inside the
region of interest.
The method can also be applied to truncated projection data by similar padding
techniques as used for filtered backprojection.

To investigate the properties of the proposed method, we computed
reconstructions using various types of prior knowledge about the reconstructed
object: box constraints on the pixel values, $\ell_1$ minimization of the
reconstruction in a wavelet basis, and $\ell_1$ minimization of the gradient of
the reconstructed image. The results show that the proposed method is able to
accurately approximate the reconstructions that would be the result of computing
the regularized iterative methods on the full object.
Compared to standard reconstruction methods such as FBP and SIRT, the proposed method is able
to significantly improve reconstruction quality by exploiting prior knowledge.

One interesting application of the method is to use it to tile reconstructions
of small subvolumes to obtain a reconstruction of the complete object. Using the
proposed method, the reconstruction of each subvolume is completely independent
of the other subvolumes. This enables parallel computation of the complete
reconstruction, resulting in a significant reduction of computation time.
The results of this paper show that the reconstruction quality of such a tiling
is comparable to the standard global regularized iterative reconstruction.
The reduction of computation time might enable the use of more advanced types
of prior knowledge that are too computationally expensive to apply globally.
Another application is to quickly estimate the parameters of a slow regularized iterative
method by estimating them in only a small subvolume.

The filter-based method of \cite{Pelt2015} on which the proposed method is based
relies on the shift-invariance of the projection operations. Therefore, it is
only applicable to parallel-beam tomography in its current form. How to apply a
similar method to other acquisition geometries is subject to
further research. It may be necessary to use additional approximations to derive
filter-based methods in other geometries, in which case exploiting prior
knowledge may actually help to reduce artifacts caused by the additional
approximations.

\section*{Acknowledgment}
This research was funded by the Netherlands Organisation for Scientific Research (NWO), project number 639.072.005.
We thank  Wolfgang Ludwig of the European Synchrotron Radiation Facility
(ESRF), Grenoble,  for providing the experimental data.
We acknowledge COST Action MP1207 for networking support.

\ifreview
\else
\fi
\bibliographystyle{IEEEtran}
\bibliography{lib}

\vfill

\end{document}